\documentclass[12pt]{article}
\usepackage{latexsym, amssymb}

\DeclareMathAlphabet\gothic{U}{euf}{m}{n}

\makeatletter
\def\eqnarray{\stepcounter{equation}\let\@currentlabel=\theequation
\global\@eqnswtrue
\tabskip\@centering\let\\=\@eqncr
$$\halign to \displaywidth\bgroup\hfil\global\@eqcnt\z@
  $\displaystyle\tabskip\z@{##}$&\global\@eqcnt\@ne
  \hfil$\displaystyle{{}##{}}$\hfil
  &\global\@eqcnt\tw@ $\displaystyle{##}$\hfil
  \tabskip\@centering&\llap{##}\tabskip\z@\cr}

\def\endeqnarray{\@@eqncr\egroup
      \global\advance\c@equation\m@ne$$\global\@ignoretrue}

\def\@yeqncr{\@ifnextchar [{\@xeqncr}{\@xeqncr[5pt]}}
\makeatother

\begin{document}

\bibliographystyle{tom}

\newtheorem{lemma}{Lemma}[section]
\newtheorem{thm}[lemma]{Theorem}
\newtheorem{cor}[lemma]{Corollary}
\newtheorem{voorb}[lemma]{Example}
\newtheorem{rem}[lemma]{Remark}
\newtheorem{prop}[lemma]{Proposition}
\newtheorem{stat}[lemma]{{\hspace{-5pt}}}
\newtheorem{obs}[lemma]{Observation}
\newtheorem{defin}[lemma]{Definition}

\newenvironment{remarkn}{\begin{rem} \rm}{\end{rem}}
\newenvironment{exam}{\begin{voorb} \rm}{\end{voorb}}
\newenvironment{defn}{\begin{defin} \rm}{\end{defin}}
\newenvironment{obsn}{\begin{obs} \rm}{\end{obs}}

\newenvironment{emphit}{\begin{itemize} }{\end{itemize}}

\newcommand{\gota}{\gothic{a}}
\newcommand{\gotb}{\gothic{b}}
\newcommand{\gotc}{\gothic{c}}
\newcommand{\gote}{\gothic{e}}
\newcommand{\gotf}{\gothic{f}}
\newcommand{\gotg}{\gothic{g}}
\newcommand{\gothh}{\gothic{h}}
\newcommand{\gotk}{\gothic{k}}
\newcommand{\gotm}{\gothic{m}}
\newcommand{\gotn}{\gothic{n}}
\newcommand{\gotp}{\gothic{p}}
\newcommand{\gotq}{\gothic{q}}
\newcommand{\gotr}{\gothic{r}}
\newcommand{\gots}{\gothic{s}}
\newcommand{\gotu}{\gothic{u}}
\newcommand{\gotv}{\gothic{v}}
\newcommand{\gotw}{\gothic{w}}
\newcommand{\gotz}{\gothic{z}}
\newcommand{\gotA}{\gothic{A}}
\newcommand{\gotB}{\gothic{B}}
\newcommand{\gotG}{\gothic{G}}
\newcommand{\gotL}{\gothic{L}}
\newcommand{\gotS}{\gothic{S}}
\newcommand{\gotT}{\gothic{T}}

\newcommand{\mn}{\marginpar{\hspace{1cm}*} }
\newcommand{\mnn}{\marginpar{\hspace{1cm}**} }

\newcommand{\mnq}{\marginpar{\hspace{1cm}???} }
\newcommand{\mnnq}{\marginpar{\hspace{1cm}**???} }

\newcounter{teller}
\renewcommand{\theteller}{\Roman{teller}}
\newenvironment{tabel}{\begin{list}%
{\rm \bf \Roman{teller}.\hfill}{\usecounter{teller} \leftmargin=1.1cm
\labelwidth=1.1cm \labelsep=0cm \parsep=0cm}
                      }{\end{list}}

\newcounter{tellerr}
\renewcommand{\thetellerr}{(\roman{tellerr})}
\newenvironment{subtabel}{\begin{list}%
{\rm  (\roman{tellerr})\hfill}{\usecounter{tellerr} \leftmargin=1.1cm
\labelwidth=1.1cm \labelsep=0cm \parsep=0cm}
                         }{\end{list}}
\newenvironment{ssubtabel}{\begin{list}%
{\rm  (\roman{tellerr})\hfill}{\usecounter{tellerr} \leftmargin=1.1cm
\labelwidth=1.1cm \labelsep=0cm \parsep=0cm \topsep=1.5mm}
                         }{\end{list}}

\newcommand{\Ni}{{\bf N}}
\newcommand{\Ri}{{\bf R}}
\newcommand{\Ci}{{\bf C}}
\newcommand{\Ti}{{\bf T}}
\newcommand{\Zi}{{\bf Z}}
\newcommand{\Fi}{{\bf F}}

\newcommand{\proof}{\mbox{\bf Proof} \hspace{5pt}} 
\newcommand{\remark}{\mbox{\bf Remark} \hspace{5pt}}
\newcommand{\ruimte}{\vskip10.0pt plus 4.0pt minus 6.0pt}

\newcommand{\simh}{{\stackrel{{\rm cap}}{\sim}}}
\newcommand{\ad}{{\mathop{\rm ad}}}
\newcommand{\Ad}{{\mathop{\rm Ad}}}
\newcommand{\Aut}{\mathop{\rm Aut}}
\newcommand{\arccot}{\mathop{\rm arccot}}
\newcommand{\capp}{{\mathop{\rm cap}}}
\newcommand{\rcapp}{{\mathop{\rm rcap}}}
\newcommand{\Capp}{{\mathop{\rm Cap}}}
\newcommand{\diam}{\mathop{\rm diam}}
\newcommand{\divv}{\mathop{\rm div}}
\newcommand{\dist}{\mathop{\rm dist}}
\newcommand{\codim}{\mathop{\rm codim}}
\newcommand{\RRe}{\mathop{\rm Re}}
\newcommand{\IIm}{\mathop{\rm Im}}
\newcommand{\Tr}{{\mathop{\rm Tr}}}
\newcommand{\Vol}{{\mathop{\rm Vol}}}
\newcommand{\card}{{\mathop{\rm card}}}
\newcommand{\supp}{\mathop{\rm supp}}
\newcommand{\sgn}{\mathop{\rm sgn}}
\newcommand{\essinf}{\mathop{\rm ess\,inf}}
\newcommand{\esssup}{\mathop{\rm ess\,sup}}
\newcommand{\Int}{\mathop{\rm Int}}
\newcommand{\Leibniz}{\mathop{\rm Leibniz}}
\newcommand{\lcm}{\mathop{\rm lcm}}
\newcommand{\loc}{{\rm loc}}

\newcommand{\mod}{\mathop{\rm mod}}
\newcommand{\spann}{\mathop{\rm span}}
\newcommand{\one}{1\hspace{-4.5pt}1}

\newcommand{\DWR}{}

\hyphenation{groups}
\hyphenation{unitary}

\newcommand{\tfrac}[2]{{\textstyle \frac{#1}{#2}}}

\newcommand{\cb}{{\cal B}}
\newcommand{\cc}{{\cal C}}
\newcommand{\cd}{{\cal D}}
\newcommand{\ce}{{\cal E}}
\newcommand{\cf}{{\cal F}}
\newcommand{\ch}{{\cal H}}
\newcommand{\ci}{{\cal I}}
\newcommand{\ck}{{\cal K}}
\newcommand{\cl}{{\cal L}}
\newcommand{\cm}{{\cal M}}
\newcommand{\cn}{{\cal N}}
\newcommand{\co}{{\cal O}}
\newcommand{\cs}{{\cal S}}
\newcommand{\ct}{{\cal T}}
\newcommand{\cx}{{\cal X}}
\newcommand{\cy}{{\cal Y}}
\newcommand{\cz}{{\cal Z}}

\newcommand{\hc}{\overline{H}}

 \thispagestyle{empty}
 
 \begin{center}
 \vspace*{-1.0cm}
 
\vspace*{1.5cm}
 
{\Large{\bf In the beginning:}  }\\[1mm]
{\Large{\bf Langlands' doctoral thesis}}  \\[4mm]
\large Derek W. Robinson$^\dag$ \\[1mm]
\large Australian National University\\[1mm]

\normalsize{March 2019}
\end{center}

\vspace{+5mm}

%
%
%



\section{Introduction}\label{S1}

Bob Langlands' mathematical research career effectively began with his 1960 PhD thesis at Yale.
It was a remarkable beginning to a remarkable career, but a beginning which largely went unnoticed.
It was remarkable as Bob wrote the  thesis, with no direct guidance,  during the first year of graduate studies; in his own words {\it ``it all happened in a hurry''}.
This did have the serendipitous outcome that during the second year of graduate work he was free to  let his interests wander in different directions.
In particular his attention focused on Selberg's work on spectral theory of Lie groups, a direction of research which led, within a few years, to the famous Langlands program. 
The thesis was also remarkable as the major part was never published in full.
The thesis consisted of two chapters.
The first chapter, approximately a third of the thesis, resolved a problem of  Hille \cite{Hil3} in the arcane area of Lie semigroups.
This material was subsequently published in the Canadian Journal of Mathematics \cite{Lan}.
The second longer and more interesting chapter only surfaced as a short
 announcement  in the Proceedings of the National Academy of Sciences \cite{Lan2}.
This announcement  is less than two pages in length and provides an excellent  illustration of   Polonius' aphorism that `Brevity is the soul of wit'.
Unfortunately the brief account failed to give an intelligible explanation of the detailed  results.
At least it was well beyond my wit and ken when I first tried to understand it in 1986.
The    second chapter 
developed a theory of general order elliptic operators affiliated with a continuous representation of a Lie group.
This theory  was then used to resolve a  problem inherent in the work of  Harish-Chandra \cite{Har} concerning the analytic structure of the representations.
The problem was very topical in the late 1950s and its solution  was a remarkable achievement for a first year graduate student
without guidance.

The Delphic nature of the Academy announcement meant that the  principal results of the thesis
passed with little attention.
Much later, in the mid 1980s,  I was interested in the integration of representations of Lie algebras and this led to an  investigation  of the regularity properties of Lie group representations.
At this point  I noticed a reference to Bob's Academy note in a paper of Roe Goodman \cite{Goo15}.
The results stated in the note were clearly of relevance to my interests at the time and I searched the literature for further evidence of the thesis. 
Frustrated by my failure to find any other published trace of the work I eventually wrote a letter to Bob asking if he still had a copy of his   thesis. 
Much to my surprise I found a copy in my mailbox a couple of weeks later.
Now, almost 60 years later, the thesis is widely available.
A photocopy of the original, which was typed by Bob's wife Charlotte in 1959, has been posted on the Princeton website {\it  http://publications.ias.edu/rpl/section/3}.  
The thesis was subsequently retyped at Yale sometime after Summer 1968  and  distributed to a  few people.
A revised and corrected version of the second typescript, prepared by myself in collaboration with Anthony Pulido, is also on the Princeton website.
An alternative more detailed presentation  of the material  in the second part of the thesis can also be found in
Chapters~I and III of  my  book `Elliptic Operators and Lie Groups' \cite{Robm}.
This latter presentation  provides an extended description  which  largely follows the reasoning of the thesis.
But both the original and my alternative version  are quite  complicated and difficult to follow.
Helgason's remark in his review of Bob's 1963 paper on automorphic forms 
{\it ``the proof involves
 many interesting ideas and techniques, which, however, do not \ldots\ldots emerge from the overly condensed exposition with the clarity they deserve''} could well apply.
Consequently the  current intention is to attempt to give a different, more pedestrian, description which explains the 
principal  difficulties and results of the thesis.
Although this description adopts the general strategy of the original and uses the same techniques
the reasoning follows a rather different route.

The story of the differential and analytic structure of representations  of Lie groups started in 1947  with  a key observation of G\aa rding \cite{Gar10}. 
He remarked that every continuous Banach space representation of a Lie group determines a representation of  the associated Lie algebra.
The image of an element of the Banach space under the action of the representation is a function on the group and 
the representatives of the Lie algebra can be viewed as  first-order differential operators with smooth coefficients acting on these functions. 
 G\aa rding  established by a regularization technique that the  corresponding subspace of  infinitely often differentiable functions is dense.
But  in 1953 Harish-Chandra \cite{Har} observed that the G\aa rding subspace was not very satisfactory  for many purposes and 
proposed restricting consideration to  the subspace of the G\aa rding space formed by the functions with convergent Taylor series.
He also introduced the terminology   well-behaved functions, or well-behaved vectors.
It was not, however, apparent that the subspace spanned by such  functions would automatically  be  dense in the representation space.
Nevertheless, Harish-Chandra  did  prove the density  for certain representations of a specific class of  groups.
His proof relied heavily on the  structure theory of Lie groups.
Subsequently, in 1958,  Cartier and Dixmier \cite{CD}  introduced the alternative terminology analytic functions
and established the density property for a wider class of groups and representations.
In particular they established that the analytic functions are dense for every unitary representation of an arbitrary Lie group.
This was the  first  indication that the density property of the analytic functions might be  a universal property, valid for general representations of arbitrary Lie groups.
This universality  was the final  conclusion, Theorem~10, of Bob's thesis.
Strangely enough the result is not explicitly stated in the Academy announcement although it is a direct inference of the final sentence.

The problem of the density of the  analytic functions  also attracted the attention of Ed Nelson who was a postdoctoral fellow 
at the Institute of Advanced Studies (IAS) in Princeton.
Nelson \cite{Nel} established the density property for general Lie groups and representations  independently of  Langlands although they were aware of each other's interests.
As Bob explained
{\it ``I do not think that Nelson ever saw my thesis. What happened $\ldots\ldots$
 is that Lennie Gross, who was then
an instructor at Yale, invited two graduate students John Frampton
and myself to drive with him to the IAS, where he was to visit friends
from his graduate student years in Chicago, among them Ed Nelson
and Paul Cohen. So it came that Nelson and I discovered a common
interest. As I have mentioned before, it was that conversation,
which I would guess went into some detail, that
provided Ed and myself with whatever we knew of each other's work.''}
The only other remaining clue to overlap of their works is a reference to a forthcoming paper of Nelson in the introduction to Bob's thesis  together with a reference to a related paper of Nelson and Stinespring \cite{NeS}.
In addition  there is a  reference to Nelson's published paper \cite{Nel}
in the Academy announcement. 
The road trip to Princeton later had an unforeseen consequence with long term implications.
As Bob wrote recently 
{\it ``I also believe that my invitation to come to
Princeton as an instructor then came to me without any further action
on my part, thus I made no application. I do not know any longer from
whom the invitation came, but it was clear that the recommendation
had come from Ed.''}
As far as I am aware Nelson never made any  specific  reference to Bob's thesis  or the announcement in the Proceedings of the Academy
in any of his publications.

The reference to Bob's Academy note  in Roe Goodmans's paper \cite{Goo15} was a result of Roe spending the academic year 1968/69 at the IAS on leave from MIT. 
By this time Bob had moved back to Yale and Nelson's interest  had moved on to Constructive Quantum Field Theory.
Fortunately Roe became aware of Bob's thesis during this period and even obtained a copy. 
 He wrote that it was {\it ``perhaps at the Institute or at the historic  old Fine Hall library''} although he subsequently admitted that 
{\it  ``perhaps it was at MIT''}.
There is now  no  trace of the thesis in the Princeton library records so the origin of the copy is obscure.
 It was possibly a copy of the second retyped version. In any case
 Roe wrote that he {\it ``read it in detail'' } and was probably the first and last reader until Bob sent me a copy in 1986.

The general strategy of Langlands and Nelson was very similar but differed significantly in detail.
Both used techniques of  semigroup theory and parabolic  partial  differential equations.
The common idea was to construct elliptic operators as polynomials in the representatives of the Lie algebra and to argue that these operators generated continuous analytic semigroups which mapped the representation space into the subspace of analytic functions.
The density property was then a direct consequence of the semigroup continuity.
The implementation of this strategy differed between the two authors. 
Nelson only considered operators in the form of Laplacians, sum of squares of the Lie algebra representatives,
but Langlands analyzed polynomials of all orders with complex coefficients.
We will not attempt to describe Nelson's arguments albeit to say that they relied in part on  probability theory and in part on results
of Eidelman  \cite{Eid} on analyticity of solutions of parabolic differential equations together with some observations 
of G\aa rding on the decrease properties of these solutions.
A complete, simplified, version of Nelson's proof was subsequently given by G\aa rding \cite{Gar11}.
Although we will not discuss Nelson's work it should be noted that it dealt with a broad range of topics involving analytic elements.
In particular it became well known  in the mathematical physics community for its results on single operator theory, the intracies involved in the addition of unbounded operators etc.
I became acquainted with it in the early 60s some 25 years before I encountered   Langlands' thesis.

The  arguments in the thesis proceeded in four major steps, the first three of which corresponded to the  theorems in the 
Academy note.
First, it was  necessary to prove that the (closure of  the) elliptic operators did indeed generate analytic semigroups.
Since Bob was considering  operators of general order this was a  very  complicated technical  problem
whose resolution was given by Theorems~7 and 8 of the thesis and stated as Theorems~1 and~2 of \cite{Lan2}.
Secondly,  as a consequence of approaching the problem in such generality, it was straightforward to deduce that 
the semigroups mapped  into the G\aa rding subspace of infinitely often differentiable functions.
Next his analysis established that the action of the semigroups was determined by a universal integrable kernel.
This was Theorem~9 of \cite{Lan} and Theorem~3 of \cite{Lan2}.
Finally it was necessary to prove that the  semigroups in fact map into the subspace of analytic functions.
This was Theorem~10 in the thesis but 
only appeared as
a passing remark at the end of the Academy paper.
It was at this last  stage that the proof depended on the theory of parabolic  partial  differential equations with analytic coefficients.
In fact Langlands cites the work of  Eidelman \cite{Eid} which was also used by Nelson.
Since Eidelman's paper  is in Russian this was a barrier in 1986 to my comprehension of the proof. 
This difficulty was compounded by  Bob's observation in his thesis  that {\it ``The facts we need from this paper are not explicitly stated as theorems and the proofs are not given in complete detail. However, since the proofs are quite complicated $\ldots\ldots$
we prefer not to perform the calculations in detail here''.}
Fortunately this final stage of the proof can be completed by a quite different argument.
The Eidelman results on the analyticity of solutions of parabolic equations were expressed in terms of complex variables
but one can instead use real variable arguments.
This approach was  given in detail in Chapter~II of \cite{Robm} and we will give a streamlined version 
 in the sequel.

Since our aim is to be as elementary as possible  we will take  a different approach to both Langlands and Nelson.
We will work in Langlands' framework with general order elliptic operators and  first give an elementary proof that the density property
of the analytic functions holds for all representations of the  multi-dimensional Euclidean group of translations.
This result is essentially a straightforward exercise in Fourier analysis.
Secondly we transfer the conclusions for the Euclidean group to all  continuous  representations of a general Lie group by using an alternative version of  the parametrix arguments developed by Langlands in combination with  some relatively simple functional analytic arguments.
In the course of our argument we establish that the elliptic operators do  generate holomorphic semigroups whose action is given by integrable kernels.
In fact we also deduce that these kernels satisfy Gaussian-style bounds.
All these conclusions are reached by variations of the arguments of Langlands' thesis supplemented by other results developed in the 1950s, the heyday of semigroup theory.
We conclude with an overview of other lines of investigation which developed from the thesis work.

\section{The Euclidean group}\label{S2}

Our discussion  of the density of analytic functions starts with an examination  of the representations of the  Euclidean group $\Ri^d$ of translations.
Since this group is commutative there are no complications of structure theory  but it nevertheless remains to establish that 
 the density property is universal for all representations of the group.
In the $\Ri^d$-case we argue that the universality follows once one has demonstrated the density of the analytic functions for the unitary representation of translations on $L_2(\Ri^d)$.
Thus the proof of the density is reduced to  understanding a relatively simple  unitary representation and then lifting, or transferring,  the result to a general representation.
A similar strategy works in the general case as we  demonstrate in the subsequent section.

We begin by recalling some well known properties of the unitary representation of $\Ri^d$ by translations on $L_2(\Ri^d)$, the Hilbert space of square integrable functions with respect to Lebesgue measure and   with norm $\|\cdot\|_2$.
The group  representation $T$  is given  by  the family of operators defined by  $(T(y)\varphi)(x)=\varphi(x-y)$ for all $x,y\in\Ri^d$ and $\varphi\in L_2(\Ri^d)$.
If $x_1,\ldots,x_d$  is a basis of $\Ri^d$ then the generators 
of the  one-parameter  subgroups $t\in\Ri\mapsto T(tx_k)$ of translations in the coordinate directions are given by 
$-\partial_k$ where $\partial_k=\partial/\partial x_k$.
The partial derivatives are the  representatives of the Lie algebra.
Next adopt the multi-index notation $x^\alpha=x_{k_1}\ldots x_{k_n}$, $\partial^\alpha=\partial_{k_1}\ldots \partial_{k_n}$ etc. where $\alpha=(k_1,\ldots,k_n)$ and the $k_j\in \{1,\ldots,d\}$.
Further denote the length of $\alpha$  by $|\alpha|=n$.

The differential structure of the $L_2$-representation is described by a well known family of Sobolev spaces.
The  subspace  $L_{2;n}(\Ri^d)$ of differential functions of order $n$ is the common domain $\bigcap_{\{\alpha:|\alpha|=n\}}D(\partial^\alpha)$ of all the $n$-th order differential operators and  the subspace of infinitely often differentiable functions is  given by $L_{2;\infty}(\Ri^d)=\bigcap_{n\geq 0}L_{2;n}(\Ri^d)$.
It follows by a standard argument with an approximate identity that $L_{2;\infty}(\Ri^d)$ is dense in $L_2(\Ri^d)$.
Finally the function  $\varphi\in L_{2;\infty}(\Ri^d)$ is defined to be an analytic function for translations  if
  $\sum_{k\geq 1}(s^k/k!)N_k(\varphi)<\infty$ for some $s>0$ and an entire analytic function if the sum converges for all $s>0$.
Here   $N_k$ is the seminorm on $L_{2;k}(\Ri^d)$ defined  by $N_k(\varphi)=\sup_{\{\alpha:|\alpha|=k\}}\|\partial^\alpha\varphi\|_2$.
This is the real analytic definition originally considered by Harish-Chandra \cite{Har}.
There is, of course, an equivalent complex analytic definition in terms of extensions of the functions $y\in\Ri^d\mapsto T(y)\varphi$
to strips in $\Ci^d$ but we will  not need to consider such extensions.
We will, however, need to norm  the subspaces $L_{2;k}(\Ri^d)$ and it is convenient  to set $\|\varphi\|_k=\sup_{\{0\leq l\leq k\}}N_l(\varphi)$.

First, following Langlands,   we introduce the $m$-th order partial differential operators $H=\sum_{\{\alpha: |\alpha|\leq m\}}c_\alpha\,(-\partial)^\alpha$ with coefficients $c_\alpha\in \Ci$ where  $m$ is an even integer.
Then $H$ is defined to be strongly elliptic if there is a $\mu>0$ such that 
\begin{equation}
\RRe\Big((-1)^{m/2}\sum_{\{\alpha: |\alpha|= m\}}c_\alpha\,\xi^\alpha\Big)\geq \mu\,|\xi|^m
\label{elrev2.0}
\end{equation}
for all $\xi\in \Ri^d$.
The largest value of  $\mu$ is called the ellipticity constant of $H$.
Thus if  $h(\xi)=\sum_{\{\alpha: |\alpha|\leq m\}}c_\alpha\,(i\xi)^\alpha$ then there are $\lambda\in\langle0,\mu]$ and $\omega\geq0$ such that $\RRe h(\xi)\geq \lambda\,|\xi|^m-\omega$
for all $\xi\in \Ri^d$.
Note that the strong ellipticity condition only involves the real part of the principal coefficients, i.e.\ those with $|\alpha|=m$.
Moreover, it is easily established that the condition is independent of the choice of coordinate basis.
A non-singular transformation of the basis does not affect the validity of the condition.
Although Langlands mainly considers strongly elliptic operators he does in part examine properties of elliptic operators
whose coefficients satisfy the weaker  condition
\begin{equation}
\Big|\sum_{\{\alpha: |\alpha|= m\}}c_\alpha\,\xi^\alpha\;\Big|\geq \mu\,|\xi|^m
\label{elrev2.00}
\end{equation}
for all $\xi\in \Ri^d$.
For brevity  we will, however, concentrate on the strongly elliptic case.

Secondly,  let $\widetilde \varphi\in L_2(\Ri^d)$ denote the Fourier transform of $\varphi\in L_2(\Ri^d)$, i.e.\
\[
\widetilde\varphi(\xi)=(2\pi)^{-d/2}\int_{\Ri^d} dx\, e^{-ix.\xi}\,\varphi(x)
\;.
\]
Then $(\widetilde{-\partial_k\varphi})(\xi)=(i\xi_k)\,\widetilde\varphi(\xi)$ and $(\widetilde{H\varphi})(\xi)=h(\xi)\,\widetilde \varphi(\xi)$, i.e.\ the elliptic differential operators act as multiplication operators on the Fourier space with multiplier $h$.
In particular $H$ is a closed operator on the subspace of functions $\varphi\in L_2(\Ri^d)$  
such that  $h\widetilde\varphi\in L_2(\Ri^d)$.
Therefore  $H$ generates a semigroup $S$ whose action is  given by
\[
(S_t\varphi)(x)=(2\pi)^{-d/2}\int_{\Ri^d} d\xi\,e^{ix.\xi}\,e^{-th(\xi)}\,\widetilde\varphi(\xi)
=(K_t*\varphi)(x)
\]
for all $t\geq 0$
where $K_t(x)=(2\pi)^{-d/2}\int_{\Ri^d} d\xi\,e^{ix.\xi}\,e^{-th(\xi)}$.
The semigroup property follows from the action by multiplication on the Fourier space and this also ensures that the kernel is a convolution semigroup, i.e.\ $K_{s+t}(x)=\int_{\Ri^d}dz\,K_s(z)\,K_t(x-z)$.

Several key properties of $S$ follow immediately from the Fourier transform definition together with Plancherel's theorem, i.e.\ the identities $\|\varphi\|_2=\|\widetilde \varphi\|_2$.
First one verifies easily that $S$ is strongly continuous, i.e.\ $\lim_{t\to0}\|(S_t-I)\varphi\|_2^2=0$.
In addition the $S_t$ satisfy the operator bounds $\|S_t\|\leq \exp(\omega t)$ with $\omega$ the constant in the  ellipticity bound on $\RRe h$.
Explicitly one has
\[
\|S_t\varphi\|_2^2=\int_{\Ri^d} d\xi\,e^{-2t\RRe h(\xi)}\,|\widetilde\varphi(\xi)|^2\leq e^{2\omega t}\,\|\varphi\|_2^2.
\]
More interestingly  one calculates that the semigroup $S$ maps $L_2(\Ri^d)$ into the subspace $L_{2;\infty}(\Ri^d)$.
For example. 
\begin{eqnarray*}
\|\partial^\alpha S_t\varphi\|_2 &\leq &\Big(\int_{\Ri^d}d\xi\,|\xi|^{2|\alpha|}e^{-2t\RRe h(\xi)}\,|\widetilde\varphi(\xi)|^2\Big)^{1/2}\\[5pt]
&\leq & t^{-k/m}\Big(\int_{\Ri^d}d\xi\,(t|\xi|^m)^{2|\alpha|/m}e^{-2(\lambda t|\xi|^m-\omega t) }\,|\widetilde\varphi(\xi)|^2\Big)^{1/2}
\leq C_k\,t^{-k/m}\,e^{\omega t}\,\|\varphi\|_2
\end{eqnarray*}
for all $\alpha$ with $|\alpha|=k$  where $C_k>0$.
The first step  uses the ellipticity bound on $h(\xi)$ and the second follows from an estimate  $(t|\xi|^m)^{2|\alpha|/m}\leq c_\lambda \exp(2\lambda t |\xi|^m)$.
Consequently $S_t$ maps $L_2(\Ri^d)$ into $L_{2;k}(\Ri^d)$ for each $k\geq1$ and hence into $L_{2;\infty}(\Ri^d)$. 
Unfortunately  this  argument does not give a good control on the constants $C_k$.
Nevertheless one can control the growth property by an  iterative argument starting from the bounds on the first derivatives $\|\partial_jS_t\varphi\|_2$.

It follows from the foregoing  that there is a $C_1>0$ such that 
\begin{equation}\label{lrev2.1}
N_1(S_t\varphi)=\sup_{j\in\{1,\ldots,d\}}\|\partial_j S_t\varphi\|_2\leq C_1\,t^{-1/m}\,e^{\omega t}\,\|\varphi\|_2
\label{elrev2.1}
\end{equation}
for all $\varphi\in L_2(\Ri^d)$.
But then
\[
N_k(S_t\varphi)\leq \sup_{j_1,\ldots,j_k}\|(\partial_{j_1} S_{t/k})\ldots (\partial_{j_k} S_{t/k})\varphi\|_2
\leq C_1^{\,k}\,(t/k)^{-k/m} e^{\omega t}\,\|\varphi\|_2
\]
for all $k\geq1$ and $\varphi\in L_2(\Ri^d)$.
Hence, by Stirling's formula, there are $a,b>0$ such that
\begin{equation}
N_k(S_t\varphi)\leq a\,b^k\,t^{-k/m}\,(k!)^{1/m}\,e^{\omega t}\,\|\varphi\|_2
\label{elrev2.10}
\end{equation}
for all $k\geq 1$, $t>0$ and $\varphi\in L_2(\Ri^d)$.
Therefore $S_t\varphi$ is an entire analytic function for translations.
But $\lim_{n\to\infty}\|S_t\varphi-\varphi\|_2=0$.
So the entire analytic functions are dense in $L_2(\Ri^d)$ and we have proved more than we set out to do.

\smallskip

After this initial skirmish with the $L_2$-representation of $\Ri^d$ we next explain how one can transfer  the density result for the analytic functions to a general continuous Banach space representation.
It is here that estimates on the semigroup kernel $K_t$ are of importance.
In addition there are two new elements entering the arguments, the continuity and the boundedness properties of the representation.
Let $\chi$ be a Banach space and $U$ a continuous representation of $\Ri^d$ by bounded operators $U(x)$, $x\in \Ri^d$, on $\chi$.
There are two types of continuity of interest, strong continuity $\|(U(x)-I)\varphi\|\to0$ as $|x|\to0$, and weak$^*$ continuity.
But a basic result of Yosida establishes that strong continuity is equivalent to weak continuity, i.e.\ equivalent to the conditions
\begin{equation}
\lim_{|x|\to0}(f, U(x)\varphi)=(f,\varphi)
\label{lrev2.2}
\end{equation}
for all $\varphi\in \chi$ and $f\in \chi*$, the dual of $\chi$.
Alternatively, if $\chi$ is the dual of a Banach space $\chi_*$, the predual of $\chi$, then $U$ is weak$^*$ continuous if $f\circ U \in \chi_*$
for all $f\in \chi_*$ and (\ref{lrev2.2}) is valid for all $\varphi\in \chi$ and $f\in \chi_*$.
Thus both types of continuity can be handled similarly.
It also follows from the group property and either form of continuity that there are $M\geq 1$ and $\rho\geq0$ such that one has
 bounds
\begin{equation}
\|U(x)\|\leq M\,e^{\rho|x|}
\label{lrev2.3}
\end{equation}
for all $x\in \Ri^d$.
Here $\|\cdot\|$ indicates the standard operator norm associated with the Banach space.

Now to emulate the earlier $L_2$-arguments for the Banach space representation it is necessary to define an analogue of the semigroup
$S$.
One direct way of doing this is by noting that on $L_2(\Ri^d)$ the semigroup satisfies
\[
(S_t\varphi)(x)=\int_{\Ri^d} dy\,K_t(y)\,\varphi(x-y)=\int_{\Ri^d} dy\,K_t(y)\,(T(y)\varphi)(x)
\;.
\]
Thus, formally at least, one has the operator representation 
\begin{equation}
S_t=\int_{\Ri^d} dy\,K_t(y)\,T(y)=T(K_t)
\label{lrev2.4}
\end{equation}
where the integral is in the weak sense.
Therefore the correct analogue of $S$ in the Banach space setting should be the semigroup
\begin{equation}
S^{(U)}_t=\int_{\Ri^d} dy\,K_t(y)\,U(y)=U(K_t)
\label{lrev2.5}
\end{equation}
where the integral is in the weak or weak$^*$ sense.
This is an observation that we will use in the subsequent discussion of general Lie groups.
But in order to make sense of  either of the relations (\ref{lrev2.4}) or (\ref{lrev2.5})  one needs control on the growth properties of 
the kernel $K_t$.
In fact the kernel and its derivatives satisfy Gaussian-type bounds.

\begin{prop}\label{plrev2.1}
There exist $b>0$ and  $\omega\geq0$, and for each  multi-index  $\alpha$ an $a_{\alpha}>0$, such  that 
\[
|(\partial^\alpha K_t)(x)|\leq a_{\alpha }\,t^{-(d+|\alpha|)/m}\,e^{\omega t}\,e^{-b(|x|^m/t)^{1/(m-1)}}
\]
for all $x\in \Ri^d$ and all $t>0$.
\end{prop}
\proof\
First consider the case that $\alpha=0$.
Then by contour integration one deduces that  
\[
|K_t(x)|=(2\pi)^{-d/2}\Big|\int_{\Ri^d} d\xi\,e^{ix.(\xi+i\eta)}\,e^{-th(\xi+i\eta)}\Big|\leq
(2\pi)^{-d/2}\int_{\Ri^d} d\xi\,e^{-x.\eta}\,\,e^{-t\RRe h(\xi+i\eta)}
\]
for all $\eta\in \Ri^d$.
But then there are $\lambda>0$ and $\sigma, \omega\geq 0$ such that 
\[
\RRe h(\xi+i\eta)\geq \lambda \,|\xi|^m-\sigma\,|\eta|^m -\omega
\;.
\]
Therefore $|K_t(x)|\leq a\, e^{\omega t}\,e^{-\eta.x}\, e^{\sigma\,t\,|\eta|^m}$ and the required bound follows by
minimizing with respect to $\eta$.

Secondly, if  $\alpha\neq0$ then the
 derivatives   introduce additional multipliers $(i\xi)^\alpha$ on the Fourier transform. 
 But for each $\varepsilon>0$ there is a $k_{\alpha, \varepsilon}>0$ such that 
\[
|(\xi+i\eta)^\alpha|\leq (t|\xi+i\eta|^m)^{|\alpha|/m}\,t^{-|\alpha|/m}\leq k_{\alpha, \varepsilon}\,t^{-|\alpha|/m}\,e^{\varepsilon t(|\xi|^m+|\eta|^m)}
\;.
\]
Then, if  $\varepsilon $ is  sufficiently small, the estimates for the derivatives follow as above. 
\hfill $\Box$

\bigskip
One also has analogous bounds on the functions $x\to x^\beta(\partial^\alpha K_t)(x)$.
\begin{cor}\label{clrev2.11}
There exist $b>0$ and  $\omega\geq0$, and for each  pair of multi-indices  $\alpha, \beta$ an $a_{\alpha,\beta}>0$, such  that 
\[
|x^\beta(\partial^\alpha K_t)(x)|\leq a_{\alpha,\beta }\,t^{-(d+|\alpha|-|\beta|)/m}\,e^{\omega t}\,e^{-b(|x|^m/t)^{1/(m-1)}}
\]
for all $x\in \Ri^d$ and all $t>0$.
\end{cor}
\proof\
The statement follows by remarking that 
for each $\varepsilon>0$ there is an $l_{\beta,\varepsilon}>0$ such that
\[
|x^\beta|\leq |x|^{|\beta|}=(|x|^m/t)^{|\beta|/m}\,t^{|\beta|/m}\leq l_{\beta,\varepsilon}\,e^{\varepsilon (|x|^m/t)^{1/(m-1)}}\,t^{|\beta|/m}
\;.
\]
Thus if $\varepsilon $ is sufficiently small  the bounds follow from those of Proposition~\ref{plrev2.1} with a slightly smaller value of $b$.
\hfill$\Box$

\bigskip
The bounds of the corollary will be applied in the sequel to differential operators whose effective order is lower than the nominal order.
If $D_n=\sum_{\{\alpha;|\alpha|\leq n\}}c_\alpha\,\partial^\alpha$ is an $n$-th order partial differential operator with smooth coefficients
supported in a compact neighbourhood of the origin and $(\partial^\beta c_\alpha)(0)=0$ for all $\beta$ with $|\beta|<(|\alpha|-k)\vee0$
then $D_n$ is defined to have effective order $k$.
Corollary~\ref{clrev2.11} then gives the bounds
\[
|(D_n K_t)(x)|\leq a_{\alpha,\beta }\,t^{-(d+k)/m}\,e^{\omega t}\,e^{-b(|x|^m/t)^{1/(m-1)}}
\]
for all $x\in \Ri^d$ and all $t>0$,
i.e. the order of the $t$ singularity is governed by the effective order $k$ rather than the real order $n$.

\smallskip

The pointwise bounds on the kernel and its derivatives allow one to deduce various weighted  bounds.
For example, the bounds of Proposition~\ref{plrev2.1} give
\begin{eqnarray*}
\sup_{x\in\Ri^d}e^{\rho|x|}|(\partial^\alpha K_t)(x)|&\leq& a_{\alpha }\,t^{-(d+|\alpha|)/m}\,e^{\omega t}\,
\sup_{x\in\Ri^d}\Big(e^{\rho|x|}e^{-b(|x|^m/t)^{1/(m-1)}}\Big)\\[5pt]
&\leq& a_{\alpha }\,t^{-(d+|\alpha|)/m}\,e^{\omega'(1+\rho^m) t}\,
\end{eqnarray*}
with $\omega'>0$.
Weighted $L_1$-bounds are also valid.
 
 \begin{prop}\label{plrev2.2}
There is an $\omega>0$ and for each multi-index $\alpha$ an $a_\alpha>0$, such that
\[
 \int_{\Ri^d}dx\,|(\partial^\alpha K_t)(x)|\,e^{\rho|x|}\leq a_\alpha \,t^{-|\alpha|/m}\, e^{\omega(1+\rho^m)t}
 \]
 for all $\rho,t>0$.
 \end{prop}
These bounds follow straightforwardly from the 
estimates of  Proposition~\ref{plrev2.1}.
 \bigskip

  The weighted estimates of Proposition~\ref{lrev2.3} combined with the continuity bounds (\ref{lrev2.3}) imply
  \[
 \Big| \int_{\Ri^d}dy\,K_t(y)\,(f,U(y)\varphi)\Big|\leq M\int_{\Ri^d}dy\,|K_t(y)|\,e^{\rho |y|}\,\|f\|\cdot\|\varphi\|
\leq a\,M\,e^{\omega(1+\rho^m)t}\,\|f\|\cdot\|\varphi\|
 \]
 for all $f\in\chi^*$, or $\chi_*$, and $\varphi\in \chi$.
 Thus $S^{(U)}_t=U(K_t)$, formally given by (\ref{lrev2.5}),  is indeed well-defined as a bounded operator on $\chi$, for each $t>0$, and one has bounds
 \[
 \|S^{(U)}_t\|\leq  a\,M\,e^{\omega(1+\rho^m)t}
 \]
 for all $t>0$.
 Then since the $K_t$ form a convolution semigroup it follows straightforwardly that  the $S^{(U)}_t$ form a continuous semigroup
 with the type of continuity dictated by the continuity of the representation, either strong  or weak$^*$.
 But one can also identify the generator of $S^{(U)}$.

 First let $X_1,\ldots, X_d$ denote the generators of the one-parameter semigroups $t\in \Ri\mapsto U(t x_k)$ in the coordinate directions
 $x_1,\ldots, x_d$, e.g.\  $X_k=\lim_{t\to0}(U(tx_k)-I)/t$.
 Then set $\chi_n=\bigcap_{\{\alpha:|\alpha|=n\}}D(X^\alpha)$ and $\chi_\infty=\bigcap_{n\geq1} \chi_n$.
 These subspaces correspond in an obvious way to the $n$-th order differentiable functions and the infinitely often differentiable functions
 respectively.
 Secondly, note that 
\begin{eqnarray*}
X_kS^{(U)}_t\varphi&=&\lim_{t\to0}\int_{\Ri^d}dy\,K_t(y)\,(U(y+tx_k)-I)\varphi/t\\[5pt]
&=&\lim_{t\to0}\int_{\Ri^d}dy\,(K_t(y-tx_k)-K_t(y))\,U(y)\varphi/t=\int_{\Ri^d}dy\,(-\partial_kK_t)(y)\,U(y)\varphi
\;.
\end{eqnarray*}
Then by iteration
\begin{equation}
X^\alpha S^{(U)}_t\varphi=U((-\partial)^\alpha K_t)\varphi=\prod_{j\in\alpha}U(-\partial_jK_{t/|\alpha|})\varphi
\label{lrev2.6}
\end{equation}
for all $\varphi\in \chi$ and all $\alpha$.
The right hand side is well-defined by another application of the bounds of Proposition~\ref{lrev2.3}.
In fact the proposition leads to bounds
\begin{equation}
\|X^\alpha S^{(U)}_t\varphi\|\leq a_\alpha\,M\,t^{-|\alpha|/m}\, e^{\omega(1+\rho^m)t}\,\|\varphi\|
\label{elrev2.2}
\end{equation}
for all $\alpha$.
Therefore one concludes that $S^{(U)}_t\chi\subseteq \chi_\infty$.
Thirdly, define the (entire) analytic functions  in $\chi$, corresponding to the representation $U$, to be the subspace of  $\varphi\in \chi_\infty$ for which $\sum_{k\geq 1}(s^k/k!)N_k(\varphi)<\infty$ for (for all) some $s>0$
where  $N_k$ is  now the seminorm on $\chi_k$ defined  by $N_k(\varphi)=\sup_{\{\alpha:|\alpha|=k\}}\|X^\alpha\varphi\|$
and the corresponding norms are again defined by $\|\varphi\|_k=\sup_{\{0\leq l\leq k\}}N_l(\varphi)$.

The principal structural result following from this discussion is the density of the entire analytic functions.

\begin{prop}\label{plrev2.3}
The operators $S^{(U)}_t=U(K_t)$ are well-defined as strong, or weak$^*$, integrals on $\chi$.
They form a strongly, or weakly$^*$, continuous holomorphic semigroup whose
generator is the strong, or weak$^*$, closure of the operator $H^{(U)}=\sum_{\{\alpha: |\alpha|\leq m\}}c_\alpha\,X^\alpha$ on $\chi_m$.
Moreover, there are $a,b>0$ and $\omega\geq 0$ such that 
\begin{equation}
N_k(S^{(U)}_t\varphi)\leq a\,b^k\,t^{-k/m}\,(k!)^{1/m}\,e^{\omega t}\,\|\varphi\|
\label{elrev2.20}
\end{equation}
for all $k\geq1$, $t>0$ and $\varphi\in\chi$.
Hence $S^{(U)}_t$ maps $\chi$ into the subspace  of entire analytic functions for $U$.
Therefore  the latter
subspace is strongly, or weakly$^*$, dense in $\chi$.
\end{prop}
\proof\
We have already discussed the definition of the $S_t$ as a continuous semigroup. 
Next the bounds on the seminorm start by iterating the bounds
\[
N_1(S_t\varphi)\leq C_1\,t^{-1/m} \,e^{\omega(1+\rho^m)t}\,\|\varphi\|
\]
 which follow for all $\varphi\in \chi$  from the estimates (\ref{elrev2.2}) with $|\alpha|=1$.
The bounds (\ref{elrev2.20}) then follow from  the factorization 
$X^\alpha S^{(U)}_t\varphi=X_{j_1}S^{(U)}_{t/n}\,\ldots X_{j_n}S^{(U)}_{t/n}\varphi$ if $\alpha=\{j_1, \ldots, j_n\}$
in direct analogy with the argument  for translations on $L_2(\Ri^d)$.
In fact the bounds (\ref{elrev2.20}) are the direct generalization of the 
 $L_2$-bounds (\ref{elrev2.10}).
The mapping property and the density property are then an immediate consequence, as before.

It remains to identify the generator of the semigroup.
But $S^{(U)}_t\chi\subseteq \chi_m=D(H_{\!(U)})$
by the mapping property.
In particular $S^{(U)}_t D(H_{\!(U)})\subseteq D(H_{\!(U)})$.
Therefore $D(H_{\!(U)})$ is a core of the semigroup generator,
 i.e.\ the generator is the closure of $H_{\!(U)}$ with respect to the norm $\|\cdot\|_m$.
\hfill$\Box$

\bigskip

The proposition  establishes the density of the analytic functions  for all the  continuous representations of $\Ri^d$ and there
are two features of the analysis which persist in the subsequent discussion of general Lie groups.
First  we have shown that the analytic properties of all continuous representations $U$ can be inferred from those of left translations
$T$  by replacing the semigroup $S_t=T(K_t)$ by the transferred operator $S^{(U)}_t=U(K_t)$.
This transference technique carries over to  the general situation.
Then properties governed by the semigroups can be   analyzed by considering  the universal semigroup  kernel on  the $L_p(\Ri^d)$-spaces.
Chapter~II of Bob  Langlands' thesis  was based on these  tactics although the presentation was somewhat different.
Our approach in the next section puts a different emphasis on the semigroup kernel and its properties.

Before proceeding to general groups we sketch a class of examples which illustrate the diversity of representations of $\Ri^d$ 
and the breadth of application of Proposition~\ref{plrev2.3}.
Let $\cl(\gotA)$ denote the space of   all bounded operators on the Hilbert space $L_2(\Ri^d)$ equipped with the usual operator norm.
This space is in fact an algebra equipped with an adjoint operator corresponding to the adjoint operation on the operators.
It is a von Neumann algebra, or a $W^*$-algebra depending on choice of terminology.
It clearly contains all bounded multiplication operators but also all bounded functions of the partial derivatives $\partial_k$.
The latter operators typically act by convolution.
Now there is an action $\tau$ of $\Ri^d$ on the algebra given by $A\in \cl(\gotA), x\in\Ri^d\mapsto \tau_x(A)=L(x)AL(x)^{-1}\in\cl(\gotA)$.
The $\tau$ form a group of $^*$-automorphisms of $\cl(\gotA)$ which is weakly$^*$ continuous.
Alternatively there are a great variety of $^*$-subalgebras which are closed in respect to the operator norm, $C^*$-algebras,
 which are invariant under $\tau$.
 On these algebras the automorphisms are weakly (strongly) continuous.
 All the foregoing results apply to this range of examples.
 My interest in the differential structure of representations of Lie groups was initially sparked by such examples and my involvement in the application of operator algebras to quantum field theory and quantum statistical mechanics (see \cite{Br1a} \cite{Br2a}).

\section{General groups} \label{S3}

The preceding analysis of representations of $\Ri^d$ gives a simple illustration of the principal results  in
 the second part  of Bob Langlands' thesis.
 In addition it provides a starting point and a strategy for the analysis
 of the continuous representations of a general Lie group $G$.
We have shown that the analytic structure of the representations of $\Ri^d$ can be inferred from the study of 
the representation of the group as left translations on the $L_p$-spaces and we
 next explain how this approach can be modified and expanded to understand the representations of the general group.
The resulting strategy  is essentially  the same  as in the  thesis but our  tactics are somewhat different.
Both the original proofs and the following arguments use three basic techniques, the exponential map, the parametrix method
and a transference technique.
We  assume the reader is conversant with the definition and the standard properties of the exponential
map and we do  not dwell on the details.
In contrast we will elaborate  on the formulation and structure of the parametrix method.
There is little to add about the transference techique as it is applied exactly as in the $\Ri^d$-case.

The next  subsection is devoted to the parametrix method as a precursor to the applications in the following two subsections.
In the first of these we  use a parametrix to establish that  each  strongly  elliptic operator affiliated with the group representation generates a holomorphic semigroup determined by an integral kernel satisfying Gaussian-type bounds.
The  proof starts from   the estimates for the Euclidean group discussed in the previous section.
Then in the following subsection we analyze the connection between the semigroups and the analytic elements
of the representation.
Throughout we use real analytic arguments and avoid complex function theory.

\subsection{The parametrix method}\label{S3.1}

Let $G$ be a connected $d$-dimensional Lie group with 
a corresponding Lie algebra $\gotg$.
The main framework of our analysis is the  representation $L$ of the group by left translations
on the space $L_2(G\,;dg)$ of functions over $G$ which are square integrable with respect
to the Haar measure $dg$.
Explicitly 
\[ 
(L(g)\varphi)(h)=\varphi(g^{-1}h)\;,
\]
for all $\varphi\in L_2(G\,;dg)$ and $g,h\in G$.
It is, however, convenient to consider in addition the representation $L$ on the associated spaces $L_p(G\,;dg)$.
The representation $L$ is usually referred to as the left regular representation although it is the direct analogue of translations
in the Euclidean case.
Note that for   $p\in[1,\infty\rangle$  the representation is strongly continuous and for  $p=\infty$
it is weakly$^*$ continuous.
All subsequent topological properties, and in particular continuity and density properties, are  understood in the strong topology
if $p\in[1,\infty\rangle$, or the weak$^*$ topology if $p=\infty$.

The first step in the analysis is to note that the Lie group is a manifold which is locally diffeomorphic to $\Ri^d$
 under the usual exponential map.
 Therefore the semigroup kernel constructed in the previous section for a given strongly elliptic operator
 on $L_2(\Ri^d)$ can be used as a local approximation for the kernel of  the corresponding elliptic operator
 on $L_2(G\,;dg)$.
 Subsequently,  starting from  the local approximation,  one can construct iteratively  a family of functions on $G$  which formally corresponds  to the semigroup kernel in the  general representation.
The iterative method  which lies at the heart of this approach is Langlands' version of Levi's parametrix method \cite{Lev}
 dating  back to 1907.
Despite its early origins the method was not well known in the 1950s and Bob no doubt learnt of it from Felix Browder's Yale lectures
which he observed to be {\it  `of a singularly instructive nature'}.
I certainly learnt of it from reading Bob's thesis 28 years later and my initial understanding was enhanced by reading Avner Friedman's book \cite{Fri} on parabolic differential equations.
The method can be formulated as a technique for solving parabolic and elliptic differential equations starting from  local approximations obtained by fixing the coefficients at a given point.
In  semigroup theory the two  approaches give approximations to the semigroup and resolvent, respectively.
These approximations are  analogous to the well known `time-dependent' and `time-independent'  methods of   perturbation theory although  one does not have a perturbation in any  conventional sense.
In the Australian vernacular it is a Claytons~$\!$\footnote{In the late 1970s when I first came to Australia there was a vigorous government campaign
against drink driving. Concurrently a non-alcoholic drink with a colour resembling whisky and appropriately bottled
was heavily promoted under the brand name Claytons `as the drink you have when you're not having a drink'.
Since then the prefix Claytons has been commonly used to indicate an ersatz product
lacking the vital ingredient.} 
perturbation theory.
We rely on the parabolic parametrix method.
It has  two remarkable features.
First it leads to a global solution starting from the initial local approximation.
Technically this arises because the terms in the series expansion are given by convolution of
terms localized in a fixed compact region.
Therefore larger distances only arise in higher order terms.
Secondly, the expansion has extremely good convergent properties.
The solution of the parabolic equation  is a function over $\Ri_+\times G$ with the variable $t\in\Ri_+$ interpretable as
the time parameter.
The series expansion is in powers  of $t$ and
 is uniformly convergent over $G$ for all $t>0$.
In many  time-dependent  problems the perturbation series are  only convergent for small times
but in the current context one obtains convergence for all times as a consequence of the Gaussian
bounds on the local Euclidean approximant.

\smallskip

The local approximation procedure starts with the exponential map.
Recall that if  $a\in\gotg$ then $\exp(a)\in G$ is defined by $\exp(a)=\gamma(1)$
where $\gamma\colon\Ri\mapsto G$ is the unique one-parameter subgroup
of $G$ whose tangent vector at the identity $e$ of $G$ is equal to $a$.
The map is a diffeomorphism of a neighbourhood of the origin in $\gotg$
to a neighbourhood of $e\in G$.
In the case of a matrix Lie group $\exp(X)$ coincides with the usual definition of the exponential
of the matrix $X$ by a power series expansion.
In the physics literature  it is commonplace 
to use the notation $e^X$ for the map since it  shares many of the basic
features of the standard exponential.
Now let $a_1,\ldots,a_d$ be a vector space basis of the Lie algebra $\gotg$ of $G$.
Then $t\mapsto \exp(-ta_j)$ is a one-parameter subgroup of $G$ and the
corresponding left translations $t\mapsto L(\exp(-ta_j))$ form a
continuous one-parameter group on each of the spaces $L_p(G\,;dg)$.
Let $A_j$ denote the generator of this group.
For example, if $G=\Ri^d$ then $A_j=-\partial_j$.
Now we consider $m$-th order operators
\[
H=\sum_{\alpha;\;|\alpha|\leq m}c_\alpha \,A^\alpha
\]
with  $c_\alpha\in\Ci$ and $m$ an even integer.
The domain of $H$ in $L_p(G\,;dg)$ is the subspace 
$L_{p;m}(G)=\bigcap_{|\alpha|\leq m}D(A^\alpha)$ of $m$-times left-differentiable
functions.
It is not difficult to establish that $L_{p;m}(G\,;dg)$ is dense in $L_p(G\,;dg)$.
Hence
$H$ is densely defined.
Then the adjoint $H^*$ of $H$ is densely defined.
Hence $H$ is closable and, for simplicity, we retain the notation $H$ for the closure.
Moreover, the subspace $L_{p;\infty}(G\,;dg)=\bigcap_{m\geq0}L_{p;m}(G\,;dg)$ of
$C^\infty$-functions is a core for each $H$.
The operators $H$ are the direct analogue of those examined in the
previous section for the Euclidean group and  the definition of strong ellipticity (\ref{elrev2.0}) 
and the definition of the ellipticity constant are  unchanged.
Again the strong ellipticity condition is independent of the choice of basis of $\gotg$ and also of the
lack of commutativity. 
If, for example,  one replaces $A_iA_j$ by $A_jA_i$  in the definition of $H$ one effectively 
introduces a modification $A_jA_i-A_iA_j$
which is linear in the $A$ by the structure relations of $\gotg$.
Therefore reordering only changes the lower order terms, those with $|\alpha|<m$, and does not 
affect the principal terms,  those with $|\alpha|=m$.
It also follows that the product $H_1H_2$ of two strongly elliptic operators with real coefficients $H_1$ and $H_2$ of orders $m_1$ and $m_2$, respectively,
is a strongly elliptic operator of order $m_1m_2$. 
The principal coefficients of the product operator are products of the principal coefficients of the component operators.

Our immediate aim is to establish that each (closed) strongly elliptic  $H$
generates a continuous semigroup $S$ on $L_2(G\,;dg)$ whose 
action  is given by   a kernel $K$ satisfying
Gaussian type bounds.
As mentioned above we approach  this problem  by constructing a family of functions $K$ by a local approximation with the kernel of the  analogous operator on the Euclidean group, corresponding formally to the semigroup kernel.
Then in the following subsection  we verify that the family does indeed have the correct properties for a semigroup
kernel and that $H$ is the generator of the semigroup.

The motivation for  the construction is the observation that the kernel $K$, if it exists,
should be a solution of the parabolic equation  
\[
(\partial_t+H)K_t=0
\]
for $t>0$ with the initial condition $K_t\to\delta$ as $t\to 0$.
Alternatively if one defines $K_t=0$ for $t\leq0$ then 
$(t,g) \mapsto K_t(g)$ from $\Ri \times G$ into $\Ci$ should be the fundamental 
solution for the heat operator $\partial_t+H$, i.e., one should have
\begin{equation}
((\partial_t+H)K_t)(g)=\delta(t) \, \delta(g)
\label{elrev3.1}
\end{equation}
for all $t\in\Ri$ and $g\in G$.
Now the parametrix method expresses $K$ as a `perturbation' expansion in  the `time' variable $t$. 

Let $\Omega\subset G$ be an open relatively compact neighbourhood of 
the identity $e\in G$
and $B_0$ an open ball in $\gotg$ centred at the origin such that
$\left.\exp\right|_{B_0} \colon B_0 \to\Omega$ is an analytic diffeomorphism.
Set $a_x=\sum_{i=1}^dx_ia_i$, for $x\in\Ri^d$, and
$B = \{ x \in \Ri^d : a_x \in B_0 \} $.
Then for
$\varphi \colon \Omega \to \Ci$ define $\hat\varphi \colon B_0 \to \Ci$ by
$\hat\varphi(x)=\varphi(\exp(a_x))$.
If $\Omega$ is small enough the image of Haar measure under this map is 
absolutely continuous with respect to Lebesgue measure.
In particular, there exists a positive
$C^\infty$-function $\sigma$ on $B$, bounded from below by a strictly positive constant,
such  that all derivatives are bounded on $B$ and such that  
\begin{equation}
\int_\Omega dg \, \varphi(g) 
= \int_{B}dx \, \sigma(x) \,{\hat \varphi}(x)  
\label{elrev3.10}
\end{equation}
for all $\varphi \in L_1(\Omega\,;dg)$.
We normalize the Haar measure $dg$ such that $\sigma(0) = 1$
and choose the modulus on $\gotg$ such that $|a_x|=|x|$ for all $x\in B$.

The key feature of the exponential map is the existence  of $C^\infty$-vector
fields $X_1,\ldots,X_{d}$ on $B$, i.e.\ first-order  partial differential operators with coefficients in $C_c^\infty(B)$, with the property
\begin{equation}
(X_k\hat\varphi)(x)
=({\widehat{A_k\varphi}})(x)
=(A_k\varphi)(\exp(a_x))
\label{elrev3.2}
\end{equation}
for all $\varphi \in C_c^\infty(\Omega)$,
where  the $A_1,\ldots,A_{d}$ are  generators of left translations.
Thus the $X_k$ are  representatives on $L_2(\Omega\,;dg)$ of the Lie algebra $\gotg$.
Moreover, 
\begin{equation}
X_k\hat\varphi=-\partial_k\hat\varphi+Y_k\hat\varphi
\label{elrev3.21}
\end{equation}
for $\varphi \in C_c^\infty(\Omega)$ where the $Y_k$ are again $C^\infty$-vector fields.
But the crucial feature is that the $Y_k$ have effective order zero as defined in Section~\ref{S2}.
Explicitly, $Y_k=\sum_{l=1}^dc_{kl}\partial_l$ with  coefficients $c_{kl}\in
C_c^\infty(B)$  which have a first-order zero at the origin.
This property is a consequence of the Baker--Campbell--Hausdorff formula.
This asserts that if $\Omega$ is sufficiently small  and $\exp(a), \exp(b)$ are both in $ \Omega$ then
 there is a $c(a,b)\in\gotg$ such that $\exp(-a)\exp(b)=\exp(c(a,b))$ and 
\[
c(a,b)=-a+b+R(a,b)
\]
where  the remainder $R(a,b)$ is a sum of  multi-commutators    of $a$ and $b$
with  each commutator containing  at least one $a$ and one $b$.
Now $X_k\hat\varphi$ is the transform of the limit as $s\to0$~of the expression
\begin{eqnarray*}
s^{-1}\Big((L(\exp(sa_k))\varphi)(\exp(a_x))-\varphi(\exp(a_x))\Big)
=s^{-1}\Big(\varphi(\exp(c(sa_k,a_x))-\varphi(\exp(a_x))\Big)
\;.
\end{eqnarray*}
The identity (\ref{elrev3.21}) follows immediately.
The term $-\partial_k$ originates from the leading term $-sa_k$ in the expression for $c(sa_k,a_x)$.
The $Y_k$, however, stem from the term linear in $s$ which occurs in the  remainder $R(sa_k,a_x)$.
This  term is of the  form $\sum^d_{l=1}c_{kl}a_l$ and, importantly,   the coefficients $c_{kl}$ have a first-order zero.
The latter property is a consequence of   the structure relations of the Lie algebra since the remainder is a sum of multi-commutators,
each of which contains at least one $a_x$.
Therefore it follows from the identity (\ref{elrev3.21}) that
\begin{equation}
{\widehat{H\varphi}}={\widehat H}_0{\hat\varphi}+{\widehat H}_1{\hat\varphi}
\label{elrev2.3}
\end{equation}
where ${\widehat H}_0=\sum_\alpha c_\alpha (-\partial)^\alpha$ is the operator with constant coefficients
corresponding to $H$ on $\Ri^d$ and ${\widehat H}_1$ is an operator of effective 
order at most $m-1$.
In particular the coefficients of  ${\widehat H}_1$ have a zero at the origin.
This local representation of $H$ on $\Ri^d$ is the starting point for constructing the semigroup corresponding to $H$ on $G$.

Let $\widetilde  K_t$ denote the kernel associated with ${\widehat H}_0$ on $\Ri^d$ but 
with $\widetilde K_t=0$ if $t\leq0$.
Further  let $\chi\in C_c^\infty(\Omega)$ with $0\leq \chi\leq 1$ and $\chi=1$ in a neighbourhood of the identity.
Then define $K^{(0)}_t$ by  setting $\widehat K^{(0)}_t={\hat\chi}\, \widetilde K_t$ on $B$.
It follows immediately from (\ref{elrev2.3}) that
\begin{eqnarray*}
((\partial_t+H)K^{(0)}_t){\widehat{\phantom{0}}}(x)&=&
((\partial_t+{\widehat H}_0)({\hat\chi}\,\widetilde K_t))(x)+
({\widehat H}_1({\hat\chi}\, \widetilde K_t))(x)\\[5pt]
&=&\delta(t) \, \delta(x)+{\widehat M}_t(x)
\end{eqnarray*}
where ${\widehat M}_t={\widehat D }\widetilde K_t$ with $\widehat D$ a partial differential operator of  the form 
$\sum_k \,{\hat\chi}_k(x)\,M^{(k)}$,
 the sum is finite, the  ${\hat\chi}_k\in C_c^\infty(B)$ and  the $M^{(k)}$ are operators  of effective 
order at most $m-1$.
Therefore the corresponding functions $K^{(0)}_t$ and $M_t$ on $G$ have compact support and satisfy the heat equation
\begin{equation}
((\partial_t+H)K^{(0)}_t)(g)=\delta(t) \, \delta(g)+M_t(g)
\label{elrev2.4}
\end{equation}
for all  $g\in G$ and $t>0$.
But it follows from the heat equations (\ref{elrev3.1}) and (\ref{elrev2.4}) that
\begin{equation}
K_t-K^{(0)}_t
=-\int^t_0ds\,K_{t-s}*M_s
\label{elrev2.80}
\end{equation}
for all $t\in \Ri$ where $*$ denotes the usual convolution  product on $G$.
Thus defining  the convolution product $\hat *$\ on $\Ri\times G$  by
\[
(\varphi\,\hat *\,\psi)_t(g)=\int_0^t ds\,(\varphi_{t-s}*\psi_s)(g)
=\int_0^t ds\int_G dh\,\varphi_{t-s}(h)\psi_{s}(h^{-1}g)
\;,
\]
one has 
\begin{equation}
K_t=K^{(0)}_t-(K\,\hat *\,M)_t
\label{elrev2.710}
\end{equation}for all $t\in\Ri$.
Then by iteration this latter equation gives 
\begin{equation}
K_t=K^{(0)}_t-(K^{(0)}\hat*\,M)_t+(K^{(0)}\hat *\,M\,\hat *\,M)_t-\ldots
\;.\label{elrev2.71}
\end{equation}
This is the parabolic parametrix expansion alluded to above. 
It represents  the solution of the parabolic equation (\ref{elrev3.1}) and 
is the analogue of the series expansion  encountered in `time-dependent'  perturbation theory.
It is the principal tool we use in the following subsection to construct the kernel $K=\{K_t\}_{t>0}$ and demonstrate that 
the family $K_t$ forms a convolution semigroup.
Although we do not  need the elliptic or `time-independent' version of the expansion we note that it can be defined  by Laplace transformation from the parabolic version.

One can in principle define  functions $L_\lambda$  for sufficiently large $\lambda$
by 
\[
L_\lambda=\int^\infty_0{\!\!dt }\,e^{-\lambda t}\,K_t
\;.
\]
Formally $L_\lambda$ is the solution of the elliptic equation 
$((\lambda I+H)L_\lambda)(g)=\delta(g)$.
Then, assuming the $L_\lambda$ are well-defined,  the parametrix identity (\ref{elrev2.710}) immediately gives the relations
\[
L_\lambda=L^{(0)}_\lambda-L_\lambda*N_\lambda
\]
where 
$L^{(0)}_\lambda$ and $N_\lambda$ are the Laplace transforms of $K^{(0)}_t$ and $M_t$,
respectively.
Now by iteration one obtains the elliptic version of the parametrix expansion
\[
L_\lambda=L^{(0)}_\lambda-L^{(0)}_\lambda*\,N_\lambda+L^{(0)}_\lambda*\,N_\lambda*\,N_\lambda-\ldots
\;.
\]
Note that it follows from the definition of the $M_t$ that $N_\lambda=DL^{(0)}_\lambda$ where $D$ is a partial differential operator of 
effective order at most $m-1$.
The unifying feature of the two versions of the parametrix method is that they both lead to inverses of partial differential operators.
The $K_t$ are inverses of the operator $(\partial_t+H)$ on $\Ri_+\times G$ and the $L_\lambda$ are inverses of the operator 
$(\lambda I+H)$ on $G$.

Next we turn to the problem of proving that the parametrix relation  (\ref{elrev2.710}) and the expansion  (\ref{elrev2.71})
are  well-defined and the $K_t$ form a convolution semigroup.

\subsection{Kernels and semigroups}\label{S3.2}

The initial step in constructing the semigroup kernel corresponding to $H$ is to prove that the expansion (\ref{elrev2.71}) determines a unique bounded integrable function $K_t$  for all $t>0$.
Subsequently we derive more detailed boundedness and smoothness property and establish that the family of functions form a convolution semigroup.

It follows from the definition of   $\widehat K^{(0)}_t$ and  $\widehat {M_t}$
that there are $a,b>0$ and $\omega\geq0$ such that
\begin{eqnarray}
|\widehat K^{(0)}_t(x)|&\leq&
a\,t^{-d/m}e^{\omega t}e^{-b(|x|^m/t)^{1/(m-1)}}
\label{elrev2.5}
\end{eqnarray}
and 
\begin{eqnarray}
|\widehat {M_t}(x)|&\leq&
a\,t^{-(d+m-1)/m}e^{\omega t}e^{-b(|x|^m/t)^{1/(m-1)}}  
\label{elrev2.6}
\end{eqnarray}
for all $x\in \Ri^d$ and all $t>0$.
These bounds follow from  Proposition~\ref{plrev2.1} and Corollary~\ref{clrev2.11} if $x\in B$ but then are obviously true 
for all  $x\in \Ri^d$ since both functions have support in $B$.
Next we convert these bounds into bounds on $K^{(0)}_t$ and $M_t$ on $G$.

One can associate with $G$ and the left-invariant Haar measure $dg$ a modulus $g\in G\mapsto |g|$
as the shortest length measured by $dg$ of the  absolutely continuous paths from $g$ to $e$.
Then the modulus is locally equivalent to the modulus on $\gotg$ by the exponential map.
In particular there is a $c>0$ such that
\[
c^{-1}|a|\leq |\exp(a)|\leq c|a|
\]
for all $a\in B_0$.
Therefore, by the choice of modulus on $\gotg$, one has 
\[
c^{-1}|x|\leq |\exp(a_x)|\leq c|x|
\]
for all $x\in B$.
In particular $|x|\geq c^{-1}|\exp(a_x)|$ for all $x\in B$.
Since $\widehat K^{(0)}_t$ and   $\widehat {M_t}$ have support in $B$ the bounds (\ref{elrev2.5}) and (\ref{elrev2.6})
immediately translate into bounds on $K_t^{(0)}$ and $M_t$.
Explicitly 
there are $a,b>0$ and $\omega\geq0$ such that
\begin{eqnarray}
|K^{(0)}_t(g)|&\leq&
a\,t^{-d/m}e^{\omega t}e^{-b\,c^{-1}(|g|^m/t)^{1/(m-1)}}
\label{elrev2.7}
\end{eqnarray}
and 
\begin{eqnarray}
|{M_t}(g)|&\leq&
a\,t^{-(d+m-1)/m}e^{\omega t}e^{-b\,c^{-1}(|g|^m/t)^{1/(m-1)}}  
\label{elrev2.8}
\end{eqnarray}
for all $g\in G$ and all $t>0$.
These estimates immediately yield the basic existence result for the $K_t$.

\begin{prop}\label{plrev3.1}
Define $K^{(n)}_t$ recursively by 
$K^{(n)}_t=-(K^{(n-1)}*\,M)_t$
where $K^{(0)}$ and $M$ are defined as above.
It follows that the series 
$\sum_{n\geq0}K^{(n)}_t$
is $L_p(G\,;dg)$-convergent to a limit $K_t\in L_p(G\,;dg)$ for all $p\in[1,\infty]$
and $t>0$. 
The function $t>0 \mapsto K_t$ is continuous  and satisfies the heat equation $(\ref{elrev3.1})$.
\end{prop}
\proof\ 
It suffices to prove that the series  is $L_1$-, and $L_\infty$-,
convergent because the $L_p$-convergence is then an immediate consequence.
The $L_1$-convergence is particularly easy because the estimates
(\ref{elrev2.7}) and (\ref{elrev2.8}) imply that
\begin{equation}
\|K^{(0)}_t\|_1\leq a\,e^{\omega t}\;\;\;\;\;{\rm and}\;\;\;\;\;
\|M_t\|_1\leq a\,t^{-(m-1)/m}e^{\omega t}
\label{elrev2.81}
\end{equation}
for suitable $a>0$ and $\omega\geq0$.
For example, one can establish analogous estimates for the $L_1(\Ri^d\,; \sigma dx)$ norms of $\widehat K^{(0)}$ and $\widehat M$ from
(\ref{elrev2.5}) and (\ref{elrev2.6}) and these translate into the $L_1(G\,;dg)$ bounds by (\ref{elrev3.10}).
But  the recursion inequalities
\begin{equation}
\|K^{(n)}_t\|_1\leq\int^t_0ds\,\|K^{(n-1)}_{t-s}\|_1\,\|M_s\|_1
\label{elrev2.811}
\end{equation}
follow from the definition of  the $K^{(n)}_t$.
Then, arguing by induction, one establishes bounds
\[
\|K^{(n)}_t\|_1\leq a\,b^n\,(t^n/n!)^{1/m}e^{\omega t}
\]
for all $n\geq0$ and all $t>0$.
In particular the series is $L_1$-convergent for all $t>0$.

The $L_\infty$-convergence is slightly more complicated.
It relies on the $L_1$-bounds (\ref{elrev2.81}) together with   the analogous $L_\infty$-bounds
\begin{equation}
\|K^{(0)}_t\|_\infty\leq a\,t^{-d/m}e^{\omega t}\;\;\;\;\;{\rm and}\;\;\;\;\;
\|M_t\|_\infty\leq a\,t^{-(m-1)/m}\,t^{-d/m}\,e^{\omega t}
\;.
\label{elrev2.820}
\end{equation}
Now one has the recursion relations
\[
\|K^{(n-1)}_{t-s}*M_s\|_\infty\leq\int_G dh|K^{(n-1)}_{t-s}(h)|\,|M_s(h^{-1}g)|\leq 
\|K^{(n-1)}_{t-s}\|_1\,\|M_s\|_\infty
\]
but these  do not immediately give useful bounds on $\|(K^{(n-1)}\hat * M)_t\|_\infty$ because the bounds (\ref{elrev2.820})
on $s>0\mapsto \|M_s\|_\infty$ are not integrable at $s=0$.
One does, however, have the alternative bounds
\[
\|K^{(n-1)}_{t-s}*M_s\|_\infty\leq \int^t_0ds\,\|K^{(n-1)}_{t-s}\|_\infty\,\sup_{g\in G}\int_G dh \,|M_s(h^{-1}g)|
\;.
\]
But now the problem is that the left-invariant measure $dh$ is not necessarily right-invariant.
Nevertheless,
\[
\int_G dh \,|M_s(hg)|=\Delta(g)^{-1}\int_G dh \, |M_s(h)|
\]
where $\Delta$ is the modular function. 
Since $M_s$ has  support in a compact $s$-independent set and $\Delta$ is locally bounded  one then concludes that
there is a $\gamma\geq1$  such that 
\[
\|K^{(n)}_t\|_\infty=\|(K^{(n-1)}\hat * M)_t\|_\infty\leq \gamma\,\int_0^t ds\,\|K^{(n-1)}_{t-s}\|_\infty\,\|M_s\|_1
\;.
\]
Then combination of these observations readily gives the recursive inequalities
\begin{equation}
\|K^{(n)}_t\|_\infty\leq \gamma \int^t_0ds\,
\Big(\|K^{(n-1)}_{t-s}\|_\infty\,\|M_s\|_1\Big)
\wedge
\Big(\|K^{(n-1)}_{t-s}\|_1\,\|M_s\|_\infty\Big)
\;.
\label{elrev2.82}
\end{equation}
These inequalities  lead to finite bounds on each  $\|K^{(n)}_t\|_\infty$   since $s\mapsto\|K^{(0)}_{t-s}\|_\infty$  is integrable at $s=0$ and $s\mapsto \|K^{(0)}_{t-s}\|_1$ is integrable at $s=t$.
 Another induction argument  indeed establishes  bounds
\[
\|K^{(n)}_t\|_\infty\leq a\,b^n\,(t^n/n!)^{1/m}\,t^{-d/m}  e^{\omega t}
\]
for suitable $a,b,\omega$, uniformly for all $t > 0$ and $n \geq 0$.
Hence one obtains uniform convergence of the series
for $K_t$.

Secondly, similar estimates allow one to verify that $t\mapsto K_t$ is continuous and satisfies the heat
equation (\ref{elrev3.1}).
For example, the heat equation is established from a term by term calculation with the expansion (\ref{elrev2.71}).
Explicitly, one has
\begin{eqnarray*}
(\partial_t+H)(K^{(0)}_t\!*\varphi)&=&\varphi+M_t*\varphi\\[8pt]
-(\partial_t+H)((K^{(0)}\,\hat *\,M)_t*\varphi)&=&-M_t*\varphi+(M\,\hat * \,M)_t*\varphi\;\;\;\mbox{ etc.}
\end{eqnarray*}
by use of the heat equation (\ref{elrev2.4}).
Addition of the terms to $n$-th order gives cancellations  leaving a single term composed of 
convolutions of $(n+1)$-factors $M_t$ with $\varphi$.
This remainder  converges to zero as $n\to\infty$ by estimates of the foregoing type.
\hfill$\Box$

\bigskip

The  proof of Proposition~\ref{plrev3.1} immediately leads to bounds
\[
\|K_t\|_1\leq a\,e^{\omega t}\;\;\;\;\;{\rm and}\;\;\;\;\;
\|K_t\|_\infty\leq a\,t^{-d/m}\,e^{\omega t}
\]
for all $t>0$.
But a slight elaboration of the proof yields much stronger results.

\begin{cor}\label{crlrev3.2}
Let $U_{\!\rho}$ denote the operator of multiplication by the function $e^{\rho|g|}$ where $\rho\geq0$.
Then there are $a>0$ and $\omega\geq0$ such that 
\begin{equation} 
\|U_{\!\rho}\, K_t\|_1\leq a\,e^{\omega(1+\rho^m) t}\;\;\;\;\;{ and}\;\;\;\;\;
\|U_{\!\rho}\, K_t\|_\infty\leq a\,t^{-d/m}\,e^{\omega(1+\rho^m) t}
\label{elrev2.9}
\end{equation}
for all   $\rho, t>0$.
Hence there are $a,b>0$ and $\omega\geq0$ such that
\begin{equation} 
|K_t(g)|\leq a\,t^{-d/m}e^{\omega t}\,e^{-b(|g|^m/t)^{1/(m-1)}}
\label{elrev2.91}
\end{equation}
for all $g\in G$ and  $t>0$.
\end{cor}
\proof\ The proof is a simple repetition of the preceding arguments  applied to the  weighted functions
$U_{\!\rho}\, K^{(n)}_t\!$.
Now
\[
|U_{\!\rho}\, K^{(n)}_t|\leq \left(|U_{\!\rho}\, K^{(n-1)}|\,\hat* \,|U_{\!\rho}\, M|\right)_t
\]
as a consequence of the recursive definition of $K^{(n)}$ and the triangle inequality for the modulus $|g|$.
Therefore the estimates (\ref{elrev2.811}) and (\ref{elrev2.82}) are  valid with $K^{(n)}$ and $M$
replaced by $U_{\!\rho}\, K^{(n)}$ and $U_{\!\rho} \,M$.
But to make recursive estimates one needs to replace the bounds (\ref{elrev2.81}) and (\ref{elrev2.820})
by bounds on $U_{\!\rho}\, K^{(0)}$ and $U_{\!\rho}\, M$.
It follows, however, from Proposition~\ref{plrev2.2}, and the discussion preceding it, that the introduction
of the weight $U_{\!\rho}$ merely introduces an additional factor $e^{\omega\rho^mt}$ to the $K^{(0)}$ and $M$ bounds,
e.g.\  one has $\|U_{\!\rho} \,K^{(0)}_t\|_1\leq a\,e^{\omega(1+\rho^m)t}$.
Therefore the induction argument now gives bounds (\ref{elrev2.9}).
But the second of these bounds also gives
\[
|K_t(g)|\leq a\,t^{-d/m}\,e^{\omega t}\,\inf_{\rho\geq0}e^{-\rho|g|}\,e^{\omega\rho^mt}
\]
which immediately yields (\ref{elrev2.91}).
\hfill$\Box$

\bigskip

Slight elaborations of these  arguments 
establish differentiability  properties of the $K_t$.
\begin{cor}\label{crlrev3.21}
The functions $K_t$ are in $L_{p;m}(G\,;dg)$  for all $p\in[1,\infty]$.
Moreover, there are $a>0$ and 
$\omega\geq0$ such that 
\[
\|U_{\!\rho}\, A^\alpha K_t\|_1\leq a\,t^{-|\alpha|/m}\,e^{\omega (1+\rho^m) t}\;\;\;\;\;{and}\;\;\;\;\;\|U_{\!\rho}\, A^\alpha K_t\|_\infty\leq a\,t^{-(d+|\alpha|)/m}\,e^{\omega(1+\rho^m) t}
\]
for all $t>0$, $\rho\geq0$ and $\alpha$ with $|\alpha|\leq m$.
\end{cor}
The $L_1$-estimates correspond to the bounds of Proposition~\ref{plrev2.2} in the Euclidean case and the $L_\infty$-estimates lead
to Gaussian-type bounds
\begin{equation}
|(A^\alpha K_t)(g)|\leq a\,t^{-(d+|\alpha|)/m}\,e^{\omega t}\,e^{-b(|g|^m/t)^{1/(m-1)}}
\label{elrev2.910}
\end{equation}
analogous to those of Proposition~\ref{plrev2.1} by the  remark at the end of the previous proof.

\medskip
\noindent{\bf Proof of Corollary~\ref{crlrev3.21}}$\;$ First note that it suffices to consider the case $\rho=0$ since the factor $U_{\!\rho}$ can be added by the argument used to
prove Corollary~\ref{crlrev3.2}.
Secondly, observe that $K^{(0)}_t\in D(A^\alpha)$ for all~$\alpha$.
Consequently  $K^{(n)}_t\in D(A^\alpha)$ for all $n$ and $\alpha$.
Moreover,
\[
A^\alpha K^{(n)}_t=-((A^\alpha K^{(n-1)})\,\hat*\,M)_t
\;.
\]
Therefore  the recursion bounds (\ref{elrev2.811}) and  (\ref{elrev2.82}) are  valid with $K^{(n)}$ replaced by $A^\alpha K^{(n)}$.
Thirdly, the $\Ri^d$-bounds of Proposition~\ref{plrev2.1} give estimates
\[
\|A^\alpha K^{(0)}_t\|_1\leq a_\alpha\,t^{-|\alpha|/m}e^{\omega t}\;\;\;\;\;{\rm and}\;\;\;\;\;
\|A^\alpha K^{(0)}_t\|_\infty \leq a_\alpha\,t^{-(d+|\alpha|)/m}e^{\omega t}
\]
for suitable $a_\alpha>0$ and $\omega\geq0$.
The derivatives $A^\alpha$ introduce the additional factors $t^{-|\alpha|/m}$.
But $s\mapsto \|K^{(0)}_{t-s}\|_1$ is integrable at $s=t$ if and only if $|\alpha|\leq m-1$.
If the latter condition is satisfied then the  sum of the sequence $A^\alpha K^{(n)}_t$ is $L_1$-convergent by the proof of Proposition~\ref{plrev3.1} to a limit $K_{\alpha;t}$ satisfying  bounds $\|  K_{\alpha;t}\|_1\leq a\,t^{-|\alpha|/m}\,e^{\omega  t}$.
Moreover,  the sum of the    $A^\alpha K^{(n)}_t$  is uniformly convergent to  $K_{\alpha;t}$  with bounds $\|K_{\alpha;t}\|_\infty\leq a\,e^{\omega t}\,t^{-(d+|\alpha|)/m}$.
Next one must  prove that $K_t$ is in $L_{\infty;m-1}(G\,;dg)$ and the sum of the series is indeed equal to $A^{\alpha} K_t$.
But left translations are weak$^*$ continuous on $L_\infty(G\,;dg)$.
Therefore   $\varphi\in D(A_j)$ if and only if $\sup_{s\in\langle0,1]}s^{-1}\|(I-L(sa_j))\varphi\|_\infty<\infty$.
Since one has $s^{-1}\|(I-L(sa_j))K^{(n)}_t\|_\infty\leq \|A_jK^{(n)}_t\|_\infty$ it immediately follows by the convergence argument that $K_t\in D(A_j)$ and $K_{\{j\};t}= A_jK_t$.
Higher order derivatives are treated by iterating this argument.

Finally one must deal with the case $|\alpha|=m$.
This is achieved by splitting the integral over $s$ into two components, the integral over $[0,t/2]$ which causes no problem, 
and the integral over $\langle t/2,t]$.
Then if $\alpha=\{k\}\cup \alpha'$ with $|\alpha'|=m-1$ one has
\[
(A^\alpha K^{(n-1)}_{t-s}*M_s)(g)=-\int_Gdh\,(A^{\alpha'}K^{(n-1)}_{t-s})(h)(A_kL(h)M_s)(g)
\;.
\]
But 
\begin{eqnarray*}
(A_kL(h)M_s)(g)&=&L(h){{d}\over{du}}(L(h^{-1})\exp(u\,a_k)L(h)M_s)(g)\Big|_{u=0}
\\[5pt]
&=&L(h){{d}\over{du}}(\exp(u\,x_k(h).a)M_s)(g)\Big|_{u=0}
\end{eqnarray*}
where $|x_k|_\infty\leq \sigma e^{\tau|h|}$ with $\sigma\geq 1$ and $\tau\geq0$.
Therefore 
\[
\|A_kL(h)M_s\|_\infty\leq \sigma\,e^{\tau|h|}\sup_{1\leq j\leq d}\|A_jM_s\|_\infty
\]
and 
\[
\int_{t/2}^tds\,\|A^\alpha K^{(n-1)}_{t-s}*M_s)\|_\infty\leq \sigma\int_{t/2}^tds\,\|U_\tau A^{\alpha'}K^{(n-1)}_{t-s}\|_\infty
\sup_{1\leq j\leq d}\|A_jM_s\|_\infty
\;.
\]
Thus  the integral is finite.
Using this decomposition gives an alternative recursion inequality which allows one to argue as before.
\hfill$\Box$

\bigskip

The foregoing argument can be  iterated to deduce that $K_t\in L_{p\,;\infty}(G\,;dg)$.  
It is a rather tedious method of deducing that the kernel is infinitely often differentiable and it
 does not give good control over the growth of the bounds as $|\alpha|$ increases.
An alternative argument can be constructed using arguments of elliptic regularity.

\smallskip

The expansion used to construct the kernel $K$ has a notable localization
feature.
The zero-order approximant $K^{(0)}$ and the remainder function $M$ are supported by $\Omega$.
Then  since  $K^{(n)}$ involves a convolution
 of $K^{(0)}$ and $n$ copies of  $M$ it must be supported by $\Omega^n$.
Thus the larger distance behaviour of the kernel is only affected by the
higher-order terms.

\smallskip

Proposition~\ref{plrev3.1} establishes that the $K_t$ satisfy the heat equation (\ref{elrev3.1})
and the expectation is they are the kernel of a holomorphic semigroup $S_t=L(K_t)$
on the $L_p(G\,;dg)$-spaces.
Our next aim is to explain this property on the Hilbert space $L_2(G\,;dg)$.
First since $\|K_t\|_1\leq a\,e^{\omega t}$ by the earlier estimates  it follows that the $S_t$
are bounded operators  on each of the $L_p$-spaces with operator norms satisfying
$\|S_t\|_{p\to p}\leq a\,e^{\omega t}$.
Moreover, one has bounds
\[
\|S_s-S_t\|_{p\to p}\leq a\,\|K_s-K_t\|_1
\]
uniformly for $s,t$ in bounded intervals of $\langle0,\infty\rangle$.
Therefore the continuity of the $K_t$ implies that the $S_t$ are uniformly continuous on bounded intervals.
This is a characteristic of holomorphic semigroups.
Moreover, it follows from the parametrix construction and the bounds on the $K^{(n)}_t$ that  
\[
\|S_t\varphi-L(K^{(0)}_t)\varphi\|_p\leq a\,t\,\|\varphi\|_p
\]
for all $\varphi\in L_p(G\,;dg)$.
Since  $L(K^{(0)}_t)\varphi\to \varphi$ as $t\to0$
it follows that  $S_t$ converges to the identity as $t\to 0$.
Hence we set $S_0=I$.
Now the main problem is to prove that the $S_t$ have the semigroup property.
But 
\[
(S_sS_t-S_{s+t})\varphi=L(K_s*K_t-K_{s+t})\varphi
\]
so it is equivalent to prove that the $K_t$ form a convolution semigroup.
This can, however, be achieved by $L_2$-arguments.

\begin{prop}\label{plrev3.2}
The operators $S_t=L(K_t)$ form a holomorphic semigroup on $L_2(G\,;dg)$
whose generator is $($the $L_2$-closure of$\,)$ $H$.
\end{prop}

The semigroup property is established by first examining 
 the real part of $H$ which is both  strongly elliptic 
and  symmetric  on  $L_2(G\,;dg)$.

\begin{lemma}\label{llrev3.3}
Each closed symmetric strongly elliptic operator $H$ on $L_2(G\,;dg)$ is  self-adjoint and lower semibounded.
\end{lemma}

The semiboundedness  is the important feature for the proof of the proposition.
Its validity was conjectured at the end Section~6 of Ed Nelson's 1959 article \cite{Nel} on analytic vectors.
Roe  Goodman raised it again in Section~3 of his 1971 article \cite{Goo15} on regularity properties of Lie group representations.
Neither author appreciated that it was a straightforward  corollary of  Langlands' parametrix arguments.
This was not observed until 1989.
It was then pointed out  in a short note  I wrote  with the late Ola Bratteli,
 Fred Goodman and Palle J\o rgensen \cite{BGJR2}.
 The main aim of the latter note was to establish a Lie group version of the classic G\aa rding inequality.
This will be discussed further in Section~\ref{S4}.

\smallskip
 
 The proof of the semiboundedness is based on resolvent arguments.

\medskip
\noindent{\bf Proof of Lemma~\ref{llrev3.3}}$\;$
Since the kernel $K_t$ corresponding to $H$ satisfies $\|K_t\|_1\leq a\,e^{\omega t}$ the Laplace transforms  
\[
L_\lambda=\int^\infty_0{\!\!dt }\,e^{-\lambda t}\,K_t
\]
are well-defined for $\lambda>\omega$ and $\|L_\lambda\|_1\leq a\,(\lambda-\omega)^{-1}$.

Moreover the operators $R_\lambda$ defined by $R_\lambda\varphi=L_\lambda*\varphi$ are bounded on  $L_2(G\,;dg)$ 
with $\|R_\lambda\varphi\|_2\leq a\,(\lambda-\omega)^{-1}\|\varphi\|_2$.
But $K_t$ satisfies the  heat equation (\ref{elrev3.1}) by Proposition~\ref{plrev3.1}.
Therefore
\begin{eqnarray*}
(\lambda I + H)R_\lambda\varphi&=&\int^\infty_0dt\,e^{-\lambda t}((\lambda I +H)K_t)*\varphi\\[5pt]
&=&\varphi+\int^\infty_0dt\,e^{-\lambda t}((\lambda -\partial_t)K_t)*\varphi=\varphi
\end{eqnarray*}
for all $\varphi\in C_c^\infty(G)$ and then by continuity 
$(\lambda I + H)R_\lambda\varphi=\varphi$
for all $\varphi\in L_2(G\,;dg)$. 
Therefore the range of $(\lambda I + H)$ is equal to  $L_2(G\,;dg)$.
This range condition together with the symmetry of $H$ implies that $H$ is  self-adjoint.

Next one has 
$(\psi, \varphi)=((\lambda I + H)R_\lambda\psi,\varphi)=
(\psi, R_\lambda^*(\lambda I + H)\varphi)$
for all $\varphi,  \psi\in D( H)$ and $\lambda>\omega$.
Therefore $\varphi=R_\lambda^*(\lambda I + H)\varphi$ and 
\[
\|\varphi\|_2=a\,(\lambda-\omega)^{-1}\|(\lambda I + H)\varphi\|_2
\]
for all $\varphi\in D( H)$.
Then it follows by spectral theory that the self-adjoint operator $ H$ is lower semibounded.
\hfill$\Box$

\bigskip

Now we are prepared to prove the proposition.

\smallskip

\noindent{\bf Proof of Proposition~\ref{plrev3.2}}$\;$
First  let $H^\dagger$ be the formal adjoint of $H$, i.e.\ the strongly elliptic  operator  with coefficients
$c^\dagger(\alpha)=(-1)^{|\alpha|}c(\alpha_*)$ where $\alpha_*=(i_n,\ldots, \alpha_1)$ when $\alpha=(i_1,\ldots, i_n)$.
Then  the real part $H_{\!R}=(H+H^\dagger)/2$ of $H$ is a symmetric  strongly elliptic operator.
Therefore Lemma~\ref{llrev3.3}  holds with $H$ replaced by  $ H_{\!R}$.
In particular  the closure of $ H_{\!R}$ is lower semibounded.
Therefore there is a $\nu\geq0$ such that 
\[
((H\varphi, \varphi)+(\varphi,  H\varphi))/2\geq-\nu\,\|\varphi\|^2_2
\]
for all $\varphi\in D( H)$.
Next observe that if $\varphi_t\in D(H)$ satisfies
the Cauchy equation
\begin{equation}
{{d}\over{dt}}\varphi_t+H\varphi_t=0
\label{elrev3.100}
\end{equation}
for all $t>0$ then
\[
{{d}\over{dt}}\|\varphi_t\|^2_2=-( H\varphi_t, \varphi_t)-(\varphi_t,  H\varphi_t)\leq
2\nu\,\|\varphi_t\|^2_2
\,.
\]
Therefore $t\mapsto e^{-\nu t}\|\varphi_t\|_2$ is a decreasing function.
Now suppose $\varphi^{(1)}_t$ and $\varphi^{(2)}_t$ both satisfy
(\ref{elrev3.100}) and $\varphi^{(1)}_t\to\varphi$, $\varphi^{(2)}_t\to\varphi$
as $t\to 0$.
Then  $\varphi^{(1)}_t-\varphi^{(2)}_t$ also satisfies the equation but  
$\varphi^{(1)}_t-\varphi^{(2)}_t\to0$ as $t\to0$.
Hence, as a consequence of the foregoing decrease property,
$\varphi^{(1)}_t=\varphi^{(2)}_t$, i.e.\  the solution of (\ref{elrev3.100}) is
uniquely determined by the initial data $\varphi=\varphi_0$. 

Now $S_tL_2(G\,;dg)\subseteq  L_{2;m}(G\,;dg)\subseteq D( H)$ by Corollary~\ref{crlrev3.21}.
Therefore 
$\varphi_t=S_{t+s}\varphi$, with $s>0$, satisfies (\ref{elrev3.100}) with
initial data $\varphi_0=S_s\varphi$ for all $\varphi\in L_2(G\,;dg)$.
Moreover, $\varphi_t=S_tS_s\varphi$ satisfies the equation with the same
initial data.
Hence 
\[
(S_{t+s}-S_tS_s)\varphi=0
\]
for all  $\varphi \in L_2(G\,;dg)$.
This establishes that $S$ is a semigroup.
But the generator $H_S$ of $S$ is an extension of $H$ on $L_{2;m}(G\,;dg)$.
Then since  $S_tL_{2;m}(G\,;dg)\subseteq  L_{2;m}(G\,;dg)$ it follows that 
 $L_{2;m}(G\,;dg)$ is a core of $H_S$.
 Therefore $H_S= H$.

Finally, the semigroup $S$ is holomorphic if and only if $S_tL_2(G\,;dg)\subseteq D( H)$ and, in addition, one has bounds
$\|H S_t\|_{2\to2}\leq a\,t^{-1}$ on the $L_2$-operator norm for all $t\in\langle0,1]$.
But both these properties are a consquence of Corollary~\ref{crlrev3.21}.
In particular
\[
 HS_t\varphi=\int_G dg\,(HK_t)*\varphi
\]
for all $\varphi\in L_2(G\,;dg)$ since $K_t\in L_{2;m}(G\,;dg)$.
Then the $L_1$-estimate of the corollary gives
\[
\|HS_t\varphi\|_2\leq a\sup_{|\alpha|\leq m}\|A^\alpha S_t\varphi\|_2\leq a'\,t^{-1}\,e^{\omega t}\,\|\varphi\|_2
\]
for all $t>0$.
Hence $S$ is holomorphic.\hfill$\Box$

\bigskip

Finally consider a general continuous  representation $U$ of $G$ on a Banach space $\chi$.
One can again associate with a fixed basis $a_1,\ldots, a_d$ of the Lie algebra $\gotg$ the generators 
$A_{U;j}$ of the one-parameter groups $t\in\Ri\mapsto U(\exp(-ta_j))$ acting on $\chi$ and 
introduce the corresponding subspaces $\chi_n=\bigcap_{\{\alpha:|\alpha|= n\}}D((A_U)^\alpha)$.
The seminorms $N_k$ are again defined on the $\chi_k$ by 
\[
N_k(\xi)=\sup_{\{\alpha:|\alpha|=k\}}\|(A_U)^\alpha\xi\|
\]
and the norms $\|\cdot\|_k$ by $\|\xi\|_k=\sup_{\{0\leq l\leq k\}}N_l(\xi)$.
Now  the strongly elliptic operator $H$ is the closure of the operator defined on $\chi_m$ in terms of the $c_\alpha$ and the $A_U$.
It again follows from the continuity of $U$ that there exist $M\geq1$ and $\rho\geq0$ such that one has bounds
\[
\|U(g)\|\leq Me^{\rho|g|}
\]
for all $g\in G$.
But the kernel $K_t$  associated with $H$ acting on $L_1(G\,;dg)$ satisfies the weighted bounds
of Corollary~\ref{crlrev3.2}.
Therefore the conclusions of Proposition~\ref{plrev3.2} for left translations $L$ of the group on $L_2(G\,;dg)$
transfer to the representation $U$ on  $\chi$.

\begin{thm}\label{tlrev3.3}
The operators $\{S^{(U)}_t\}_{t>0}$ defined by
\[
S^{(U)}_t=U(K_t)=\int_Gdg\,K_t(g) \, U(g)
\]
form a holomorphic  semigroup 
$S^{(U)}$ on  the Banach space $\chi$ such that $S^{(U)}_t\chi\subseteq \chi_\infty$ for all $t>0$.
The generator of the semigroup is the closure of  $H$ on $\chi_m$.
\end{thm}
\proof\ 
The $S^{(U)}$ are well defined as bounded operators on $\chi$ by the bounds of Corollary~\ref{crlrev3.2}.
Then since $S$ is a semigroup for the left regular representation $L$ of $G$ on $L_2(G\,;dg)$, by Proposition~\ref{plrev3.2}, 
the $K_t$ must form  a  convolution semigroup, i.e.\ $K_{s+t}=K_s*K_t$ for all $s,t>0$.
Therefore $S^{(U)}$ is a continuous semigroup on $\chi$.
The remaining statements follow from the properties of the kernel $K$ once one notes that 
\[
U(\exp(sa))S^{(U)}_t\xi=\int_Gdg\,K_t(g)U(\exp(sa)g)\xi=\int_Gdg\,K_t(\exp(-sa)g))U(g)\xi
\]
for $a\in\gotg$, all small $s$ and all $\xi\in \chi$.
Thus $A_US^{(U)}_t\xi=U(AK_t)\xi$ for all $\xi\in\chi$  where the $A_U$ are generators
in the representation $U$ and the $A$ are the left regular representatives.
Then by iteration $(A_U)^{\alpha} S^{(U)}_t=U(A^{\alpha} K_t)\xi$ for all $\alpha$.
Note that 
$K_t$ is $m$-times left differentiable, again by  Corollary~\ref{crlrev3.21}, but in fact it is infinitely
often differentiable by the remark following the proof of the corollary.
Hence  $S^{(U)}_t\chi\subseteq \chi_m\subseteq D(H)$.
Since $\chi_m$ is $S^{(U)}$-invariant it is a core of the semigroup generator.
Thus the generator must be $H$. 
Finally the holomorphy of $S^{(U)}$ follows by the reasoning in the proof of 
 Proposition~\ref{plrev3.2}.
\hfill$\Box$

\bigskip

One also has the following extension of Lemma~\ref{llrev3.3}.

\begin{cor}\label{clrev3.31}
Assume $\chi$ is a Hilbert space,  $U$ a unitary representation and the strongly elliptic operator $H$ corresponding to $U$ is symmetric.
Then  $H$ is self-adjoint and lower semibounded with a  bound which  is independent of the choice of unitary representation.
\end{cor}
\proof\
Let $K_t$ be the kernel corresponding to $H$ and set $S^{(U)}_t=U(K_t)$.
It follows by the proof of Lemma~\ref{llrev3.3} that the generator of $S^{(U)}$ is equal to $H$ and is self-adjoint.
But  $\|U(g)\|=1$ for all $g\in G$ by unitarity where $\|\cdot\|$ the Hilbert space operator norm.
Then $\|S^{(U)}_t\|\leq \|K_t\|_1\leq a\,e^{\omega t}$ for all $t>0$.
Hence the generator  of $S^{(U)}$ is self-adjoint and lower semibounded by spectral theory.
The bound does not depend on $U$.
It is equal to the bound on $L_2(G\,;dg)$.
\hfill$\Box$

\bigskip

Another  consequence of Proposition~\ref{plrev3.2} is that the parametrix function $L_\lambda$ is the kernel of the 
resolvent of $ H$  on $L_2(G\,;dg)$.
Explicitly
\[
L_\lambda*\varphi=\int^\infty_0{\!\!dt }\,e^{-\lambda t}K_t*\varphi=\int^\infty_0{\!\!dt } \,S_t\varphi=(\lambda I+ H)^{-1}\varphi
\]
for all $\varphi\in L_2(G\,;dg)$ and $\lambda>\omega$.
Then it follows by a similar calculation that $U(L_\lambda)$ is the kernel of the generator  of the semigroup $S^{(U)}$ on $\chi$.
Thus the functions $K_t$ and $L_\lambda$ determine the action  of the semigroup and  resolvent  for all continuous representations of the group.

\smallskip

Theorem~\ref{tlrev3.3} encompasses two of the principal results of Langlands' 1960 thesis,  Theorems~8 and 9. 
Elaborations of these theorems are given by Theorem~I.5.1 and Theorem~III.2.1 in my 1991 book \cite{Robm}.
The  common  features of these earlier descriptions and the  current presentation are the exponential map and the parametrix technique but the proofs differ in some features.
For the interested reader we make a few remarks on the relations between the three versions.

In our current notation Theorem~8 of the thesis established that the closure of $H$ in the continuous representation $U$
generates a holomorphic semigroup $S^{(U)}$.
Then Theorem~9 showed  that the action of $S^{(U)}$ is given by a family of measures which form  a convolution semigroup.
Subsequently, after the statement of Theorem~10, it was shown that these measures are absolutely continuous with respect to Haar
measure, i.e.\ the  $S^{(U)}$  has a kernel $K^{(U)}$ in the current terminology.
Theorem~10 derives the basic complex analytic properties  of the semigroup.
We have bypassed Theorem~7 of the thesis which gives an identification
of the adjoint of $H$ in the representation $U$.
In fact the identification is valid for the  wider class of elliptic operators with coefficients satisfying (\ref{elrev2.00}).
Strong ellipticity is not necessary.
The proof, however, requires the parametrix expansion for the resolvent and a version of  elliptic regularity
for operators with continuous coefficients. 
This approach also leads to the alternative proof that $K_t\in L_{p;\infty}(G\,;dg)$ alluded to after the proof of Corollary~\ref{crlrev3.21}.
The original proof of the generator property in Theorem~8  of the thesis utilized Theorem~7 and, of course,  required strong ellipticity.
We have avoided these arguments for the sake of brevity.
Bob  cites Theorem~12.8.1 of Hille and Phillips book \cite{HP} in the derivation of the generator property although this appears to require some  poetic license. 
Further details can be found in Sections~I.4 and I.5 of \cite{Robm}.
We note that the identification given by the missing Theorem~7 is a  direct corollary of Theorem~\ref{tlrev3.3} 
if $H$ is strongly elliptic. 

The main distinction between the developments of the thesis and the current exposition is  largely in emphasis.
The primary focus of the thesis is  the semigroup $S^{(U)}$ and the kernel $K^{(U)}$ is a secondary artefact.
In the preceding analysis, however,  the kernel plays the principal role and the semigroup is constructed from the kernel.
The current presentation is based on three papers I coauthored \cite{BGJR2} \cite{Rob14} \cite{ER14} some 30--35 years
 after Bob wrote his thesis.
 Nevertheless, the additional arguments  are all from the late 1950s.
For example, the uniqueness criteria for the Cauchy problem are discussed at length in Chapter~III of Hille and Phillips 1957 book \cite{HP}
and the characterization of holomorphy of $S$ by a norm estimate on $H S_t$ was a 1958 result of Yosida \cite{Yos1}.
Finally the characterization of cores of semigroup generators  by semigroup invariance is essentially contained in Nelson's 1959 paper \cite{Nel}.

We will return to the discussion of other properties that can be deduced by variations of these arguments in Section~\ref{S4}.
But next we turn to the explanation of the last result, Theorem~10, of Langlands' thesis.

\subsection{Differential and analytic structure}\label{S3.3}

Theorem~10 of the thesis established  the result highlighted in the introduction,  the density of the analytic elements $\chi_a(U)$ of  a general continuous representation $U$  of a Lie group.
Now
\[
\chi_a(U)=\{\xi\in \chi_\infty:\sum_{k\geq 1}(s^k/k!)N_k(\xi)<\infty\ \mbox{ for some } s>0\}
\]
where the seminorms are defined in terms of the representatives $A_U$ of the Lie algebra.
The strategy of the thesis was to prove that $\chi_a(U)$ contains the dense subspace $\chi_a(H)$ of analytic elements of  each strongly elliptic operator $H$, i.e.\ the subspace
\[
\chi_a(H)=\{\xi\in \chi_\infty:\sum_{k\geq 1}(s^k/k!)\|H^k\xi\|<\infty\ \mbox{ for some } s>0\}
\;.
\]
The density of the latter subspace follows because  $H$ generates the holomorphic semigroup $S^{(U)}$, by Theorem~\ref{tlrev3.3}.
In particular it follows that 
\[
\chi_a(H)=\bigcap_{t>0} S^{(U)}_t\chi
\;.
\]
The density of $\chi_a(H)$  is then an immediate consequence since $S^{(U)}_t\to I$ as $t\to0$.
We give a different proof of this result based on arguments of  a 1988 paper \cite{BGJR1} whose principal focus was
the integrability of representations of Lie algebras.
This is the original problem which attracted my attention to the Lie group theory  and led me to search out Bob's thesis.
The proof uses real analytic arguments which replace the application of Eidelman's results on the analyticity of the solutions
of parabolic equations cited in the thesis. 
This application is by Bob's own admission, quoted in the introduction, somewhat nebulous.
The alternative proof which we next outline is given in greater detail in Chapter~II of \cite{Robm}.
It is  a basically a straightforward  exercise in functional analysis.
 
\begin{thm}\label{tlrev3.21}
There exist $a,b>0$ and $\omega\geq0$ such that
\begin{equation}
\|S^{(U)}_t\xi\|_n\leq a\,b^n\, n!\,t^{-n/m}\,e^{\omega t}\,\|\xi\|
\label{elrev3.22}
\end{equation}
for all $\xi\in\chi$,  $n\in\{0,1,2\ldots\}$ and $t>0$.
Therefore  $S^{(U)}_t\chi\subseteq \chi_a(U)$ for all $t>0$.
Consequently $\chi_a(H)\subseteq  \chi_a(U)$  and $\chi_a(U)$ is dense in $\chi$.
\end{thm}
\proof\
Throughout the proof we omit the index $U$.
This should cause no confusion as all calculations are within the representation.
We also use the convention that $a, b, \omega$ etc. are $n$-independent constants whose values might vary line by line.
Any $n$-dependence is explicitly noted by  suffices.

First observe that  $S^{(U)}_t\chi\subseteq \chi_\infty$ for all $t>0$ by Theorem~\ref{tlrev3.3}.
Therefore it is sufficient to derive the estimates (\ref{elrev3.22}) for all $\xi\in \chi_\infty$.
In particular if  $\lambda\in\langle0,1\rangle$ then  $\|S_t\xi\|_n=\|S_{(1-\lambda)t}(S_{\lambda t}\xi)\|_n$ with  $S_{\lambda t}\xi\in\chi_\infty$ satisfying bounds $\|S_{\lambda t}\xi\|\leq M\,e^{\omega\lambda t}\|\xi\|$.
Therefore the bounds (\ref{elrev3.22}) for general $\xi\in\chi$ follow from those for $\xi\in\chi_\infty$, albeit with increased values of $a,b$ and $\omega$.

Secondly, it suffices to establish (\ref{elrev3.22}) for all $t\in\langle0,t_0]$ for some $t_0>0$.
This follows since   $\|S_t\xi\|_n\leq (a\,b^n\,n!)\|S_{t-t_0}\xi\|$  for $t>t_0$  as a consequence of  the bounds for $t\leq t_0$.
But $\|S_{t-t_0}\xi\|\leq M\,e^{\omega(t-t_0)}\|\xi\|\leq M\,t^{-n/m}\,e^{\omega t}\,\|\xi\|$ for $t>t_0$.
Therefore one obtains the bounds (\ref{elrev3.22}) for all $t>0$, again with increased values of $a,b$ and $\omega$.
Hence the proof is now reduced to considering $\xi\in\chi_\infty$ and small $t>0$.
The starting point is the following observation for small $n$.

\begin{lemma}\label{llrev3.22}
There is an $a>0$ such that  $\|S_t\xi\|_n\leq a\,t^{-n/m}\,\|\xi\|$ for all $\xi\in \chi$, $t\in\langle0,1]$ and $n=0,1,\ldots, m-1$.
\end{lemma}
\proof\
It follows from the $L_1$-estimates of Corollary~\ref{crlrev3.21} that 
\begin{eqnarray*}
\|S_t\xi\|_n=\sup_{\{\alpha:|\alpha|\leq n\}}\|A^\alpha S_t\xi\|&=&\sup_{\{\alpha:|\alpha|\leq n\}}\|U(A^\alpha K_t)\xi\|\\[5pt]
&\leq& M\,\sup_{\{\alpha:|\alpha|\leq n\}}\|U_\rho A^\alpha K_t\|_1\|\xi\|\leq a\,t^{-n/m}\,\|\xi\|
\end{eqnarray*}
by the $(M,\rho)$-continuity estimates for the representation.
\hfill$\Box$

\bigskip

If, however,  $\xi\in \chi_\infty$ then  $\|S_t\xi\|_n$ is not singular at $t=0$.

\begin{lemma}\label{llrev3.221}
There is a $b>0$ such that  $\|S_t\xi\|_n\leq a\,(\|H\xi\|+\|\xi\|)$ for all $\xi\in \chi_\infty$, $t\in\langle0,1]$   and $n=0,1,\ldots, m-1$.
\end{lemma}
\proof\
Since  one has
\begin{eqnarray*}
\|A^\alpha(\lambda I+H)^{-1}\xi\|&\leq& \int^\infty_0ds\,e^{-\lambda s}\,\|A^\alpha S_s\xi\|\\[5pt]
&\leq & \int^\infty_0ds\,e^{-\lambda s}\,\|U(A^\alpha K_s)\xi\|
\leq a\int^\infty_0ds\, e^{-(\lambda-\omega)s}s^{-n/m}\,\|\xi\|
\end{eqnarray*}
 the statement of the lemma follows by replacing $\xi$ by $(\lambda I+H)S_t\xi$
and choosing $\lambda$ sufficiently large.
\hfill$\Box$

\bigskip

The rest of the proof consists of `bootstrapping' the small $n$-estimates of Lemma~\ref{llrev3.22}
  into universal estimates with the correct quantitative behaviour 
for small $t>0$.
This is achieved by recursive arguments.

Assume one has  bounds 
\begin{equation}
\|S_t\xi\|_n\leq c_n\,t^{-n/m}\,\|\xi\|
\label{elrev3.220}
\end{equation}
with $c_n>0$
for all $\xi\in \chi_\infty$, $t\in\langle0,1]$ and all $n\leq  k(m-1)$  for some $k\geq 1$.
 It follows from the lemma this assumption is indeed valid for $k=1$, with $c_1=\ldots=c_{m-1}=a$.
Next we argue that similar bounds are valid  for all $ n\leq (k+1)(m-1)$ and we also obtain  estimates on the corresponding $c_n$.
A crucial part of the argument is the observation that the commutators
$(\ad\, A^\alpha)(A^\beta)=A^\alpha A^\beta-A^\beta A^\alpha$ have the property
\[
\|(\ad\, A^\alpha)(A^\beta)\xi\|\leq a_{\alpha, \beta}\|\xi\|_{\alpha+\beta-1}
\]
for all $\alpha,\beta$ and all $\xi\in\chi_\infty$ as a result of the structure relations of $\gotg$.
In particular the commutator gives a unit reduction to the apparent operator order.

Let $A^\alpha=A^{\alpha_0}A^{\alpha_1}$ with $|\alpha|=n\geq m$ and $|\alpha_0|=m-1$ so $|\alpha_1|=n-m+1$.
Then
\begin{eqnarray*}
A^\alpha S_t\xi&=& A^{\alpha_0}S_s(A^{\alpha_1}S_{t-s}\xi)+A^{\alpha_0}(\ad\, A^{\alpha_1})(S_s)S_{t-s}\xi\\[5pt]
&=& A^{\alpha_0}S_s(A^{\alpha_1}S_{t-s}\xi)-\int^s_0du\,A^{\alpha_0}S_u(\ad\, A^{\alpha_1})(H)S_{t-u}\xi
\end{eqnarray*}
for all $s\in\langle0,t\rangle$ and $\xi\in \chi_\infty$.
Since $H$ is a polynomial in $A^\beta$ with $|\beta|\leq m$ one then deduces from (\ref{elrev3.220}) that 
\begin{eqnarray*}
\| S_t\xi\|_n&\leq & a\,c_{n-m+1}\,s^{-(m-1)/m}(t-s)^{-(n-m+1)}\,\|\xi\|\nonumber \\[5pt]
&&\hspace{3cm}{}+b\,n\int^s_0dr\,r^{-(m-1)/m}\,\|S_{t-r}\xi\|_n
\end{eqnarray*}
where $a, b>0 $ are independent of $n$.
Now it remains to solve these inequalities.

First set $s= \mu t$ with $\mu\in\langle0,1\rangle$ and $r=ut$.
Then
\begin{eqnarray*}
\|S_t\xi\|_n&\leq & a\,c_{n-m+1}\,t^{-n/m}\,\mu^{-(m-1)/m}(1-\mu)^{-(n-m+1)}\,\|\xi\|\nonumber \\[5pt]
&&\hspace{3cm}{}+b\,n\,t^{1/m}\int^\mu_0du\,u^{-(m-1)/m}\,\|S_{t(1-u)}\xi\|_n
\;.
\end{eqnarray*}
Second, choose $\mu=n^{-m}$. 
Then $\mu^{-(m-1)/m}=n^{m-1}$ and $(1-\mu)^{-(n-m+1)}\leq b_m$ where $b_m=\sup_{n\geq m}(1-n^{-m})^{-n/m}<\infty$.
Therefore 
\begin{eqnarray*}
\|S_t\xi\|_n\leq  a\,n^{m-1}\,c_{n-m+1}\,t^{-n/m}\|\xi\|
+b\,n\,t^{1/m}\int^{n^{-m}}_0du\,u^{-(m-1)/m}\,\|S_{t(1-u)}\xi\|_n
\;.
\end{eqnarray*}
But iteration of this inequality $p$-times gives
\begin{eqnarray*}
\|S_t\xi\|_n&\leq & a\,n^{m-1}\,c_{n-m+1}\,t^{-n/m}\,\|\xi\|\;\sum^{p-1}_{l=0}(c\,t^{1/m})^l+R_{p,n}(t)
\end{eqnarray*}
with $c=bm$ independent of $n$ and 
\[
R_{p,n}(t)=(b\,n\,t^{1/m})^p\int^{n^{\!-m}}_0\!\!\!\!du_1\,u_1^{-(m-1)/m}\ldots \int^{n^{\!-m}}_0\!\!\!\!du_p\,u_p^{-(m-1)/m}\,\|S_{t(1-u_1)\ldots(1-u_p)}\xi\|_n
\;.
\]
But  if $t<t_0$ with $c\,t_0<1$ then in the limit  $p\to\infty$  one has 
\[
\|S_t\xi\|_n\leq  a\,n^{m-1}\,c_{n-m+1}\,t^{-n/m}\,\|\xi\|+\limsup_{p\to\infty}R_{p,n}(t)
\]
for all $t<t_0$ and $\xi\in\chi_\infty$.
Next we argue that the limit of the remainder term is zero.

We may assume $t\leq t_0<1 $  and $n\geq 1$.
Therefore one immediately has bounds
\[
R_{p,n}(t)\leq (b\,t_0)^p\,\bigg(n\int^{n^{\!-m}}_0\!\!\!\!du\,u^{-(m-1)/m}\bigg)^p\sup_{t\in\langle0,1]}\|S_t\xi\|_n=
(c\,t_0)^p\sup_{t\in\langle0,1]}\|S_t\xi\|_n
\;.
\]
Now $n\leq (k+1)(m-1)<(k+1)m$.
Then, by a relatively straightforward extension of Lemma~\ref{llrev3.221}, one obtains bounds $\sup_{t\in\langle0,1]}\|S_t\xi\|_n\leq a_k\,(\|H^{k+1}\xi\|+\|\xi\|)$, e.g.\ if the principal coefficients of $H$ are real then $H^{k+1}$ is strongly elliptic and the lemma is valid with $H$ replaced by $H^{k+1}$ and $m$ replaced by $(k+1)m$.
Hence 
\[
R_{p,n}(t)\leq a_k\,(c\,t_0)^p\,(\|H^{k+1}\xi\|+\|\xi\|)
\]
for all $\xi\in \chi_\infty$.
Thus $R_{p,n}(t)\to0$ as $p\to\infty$ for all $t\in\langle0,t_0\rangle$ if $c\,t_0<1$.
Combining these observations one concludes that 
\[
\|S_t\xi\|_n\leq  a\,n^{m-1}\,c_{n-m+1}\,t^{-n/m}\,\|\xi\|
\]
for all $t<t_0$, with $t_0$ sufficient small,  and all $\xi\in\chi_\infty$.
Moreover, by induction, these estimates are valid for all $n\geq m-1$.

Next if  $k(m-1)\geq n>(k-1)(m-1)$ then by iteration
\[
c_n\leq a^{k-1}\,\Big(\prod^{k-2}_{l=0} (n-l(m-1))^{m-1}\Big)\, c_{n-(k-1)(m-1)}\leq  a^{k-1}\,n^n\, c_{m-1}
\]
because the product has $k-1$-factors each bounded by $n^{m-1}$.
But  $k-1<n/(m-1)$ and $n^n\leq e^n\,n!$.
Thus the bounds take the form $a\,b^n\,n!$ for all $n\geq 1$.

Finally, since $N_n(\xi)\leq \|\xi\|_n$ for all $\xi\in\chi_\infty$ one has
\[
\sum_{k\geq1}(s^k/k!)\,N_k(S_t\xi)\leq \sum_{k\geq1}(s^k/k!)\,\|S_t\xi\|_k\leq a\sum_{k\geq1}(bst^{-1/m})^ke^{\omega t}\|\xi\|
\]
for all $\xi\in \chi$ and $bs<t^{1/m}$. 
Therefore $S_t\chi\subseteq \chi_a(U)$ for all $t>0$.
\hfill$\Box$

\bigskip

Theorem~\ref{tlrev3.21} achieves the aim set out in the introduction, it establishes the density of the analytic elements
for all continuous representations of a general Lie group.
This property was the final conclusion, Theorem~10,  of Langlands' 1960 thesis.
It is also the final conclusion of our explanation of the results of the thesis.
But we are not finished. 
To conclude we give a summary of consequences of the thesis results which have been subsequently established.
Most of these results date from the Gaussian revolution in semigroup theory which started slowly in 1967  with Aronson's paper
\cite{Aro} on bounds on solutions of parabolic equations and which peaked in the 1980s.

\section{Consequences}\label{S4}

The early results of Langlands on the differential and analytic structure of Lie groups have developed in two different, but related, frameworks.
First, there has been considerable progress in the framework of strongly elliptic operators described above.
Secondly, the theory has been generalized to a broader class of subelliptic operators.
The latter operators are defined as polynomials in an algebraic basis of the Lie algebra $\gotg$, i.e., a linearly independent subset of $a_1,\ldots a_{d_1}\in \gotg$  whose Lie algebra spans $\gotg$.
We will briefly describe some of the key features of the developments in the strongly elliptic setting and then comment on the more complicated subelliptic situation.
We have seen in the foregoing that despite appearances the strongly elliptic theory remains largely commutative. 
The subelliptic theory, in contrast,  contains a genuine noncommutative element.

The first topic of our discussion is  extensions of the classical G\aa rding inequality to unitary representations of Lie groups.
These inequalities characterize the notion of  strong ellipticity   and  provide a basis for the definition of  subellipticity for  operators
of general order.

\begin{prop}\label{plrev4.1}
Let $U$ be a unitary representation of $G$ on a Hilbert space $\chi$ and $H$ a strongly elliptic operator with ellipticity constant $\mu$.
Then for each $\lambda\in \langle0,\mu\rangle$ there is a representation independent $\nu\geq0$ such that 
\begin{equation}
\RRe(\varphi, H\varphi)\geq \lambda\, N_{m/2}(\varphi)^2-\nu\,\|\varphi\|^2
\label{elrev4.1}
\end{equation}
for all $\varphi\in \chi_m$.
\end{prop}
\proof\
Since $H$ is strongly elliptic  with ellipticity constant $\mu$ it follows that  the real part of  $H-\lambda\,A^{\alpha_*}A^\alpha$, with $|\alpha|=m/2$,  is a symmetric  strongly elliptic  operator with ellipticity constant $\mu-\lambda$.
Therefore  it follows from Corollary~\ref{clrev3.31} that $\RRe H-\lambda\,A^{\alpha_*}A^\alpha$ is essentially self-adjoint and lower semibounded by $-\nu_\alpha I$ where
 $\nu_\alpha=\inf_{t>0}\log\|K^{(\alpha)}_t\|_1$ and $K^{(\alpha)}_t$ is the kernel corresponding to the operator $\RRe H-\lambda\,A^{\alpha_*}A^\alpha$.
Hence
\[
\RRe(\varphi, H\varphi)\geq \lambda\,\|A^\alpha\varphi\|^2-\nu_\alpha\,\|\varphi\|^2
\]
for all $\varphi\in \chi_m$.
 But $\nu_\alpha$ is  independent of the particular unitary representation.
Therefore (\ref{elrev4.1}) follows by taking the supremum over the $\alpha$.
\hfill$\Box$

\bigskip

The inequality (\ref{elrev4.1})  is a Lie group version of the classic G\aa rding inequality for strongly elliptic divergence form operators with bounded continuous coefficients on $L_2(\Ri^d)$.
But there are other possible formulations.
For example,
\begin{equation}
\RRe(\varphi, H\varphi)\geq \lambda\, (\varphi,
\Delta^{m/2}\varphi)-\nu\,\|\varphi\|^2
\label{elrev4.0}
\end{equation}
for all $\varphi\in \chi_m$ with $\Delta$  the Laplacian corresponding to the 
basis $a_i$ in the representation $U$.
The proof follows as before since $H-\lambda\, \Delta^{m/2}$  is strongly elliptic for all $\lambda\in \langle0,\mu\rangle$.
 
 \smallskip

Proposition~\ref{plrev4.1} establishes that strong ellipticity of $H$ implies the G\aa rding inequalities~(\ref{elrev4.1})
for any unitary representation of $G$.
Moreover, since the the strong ellipticity condition (\ref{elrev2.0}) is just a restriction on the coefficients $c_\alpha$
of the operator  it also implies the G\aa rding inequalities for the unitary representation of $\Ri^d$ by left translations on $L_2(\Ri^d)$.
Conversely assume that  (\ref{elrev4.1}) is valid without the strong ellipticity 
 restriction (\ref{elrev2.0}).
Then  choose  $\varphi\in C^\infty_c(\Omega)$ 
with ${\hat\varphi}(x)=e^{i\eta \cdot x}\chi(x)$ where  $\chi$ is a $C^\infty$-function
 supported in a ball of radius $r$ centred at the origin and we have again 
used  the exponential map and the notation of Subsection~\ref{S3.1}.
Evaluating   (\ref{elrev4.1})  with 
$\varphi$ for large $\eta$ and small $r$ then yields the strong ellipticity
condition. 
These arguments are summarized by the following.

\begin{cor}\label{clrev4.10}
The strong ellipticity condition $(\ref{elrev2.0})$ is equivalent to  the G\aa rding
inequality $(\ref{elrev4.1})$
in a unitary representation of $G$, or
equivalent to $(\ref{elrev4.1})$ for left translations on $L_2(\Ri^d)$.
\end{cor}

A similar conclusion follows with the G\aa rding inequality (\ref{elrev4.1}) replaced by the alternative formulation (\ref{elrev4.0}).
In particular both forms of the G\aa rding inequality are equivalent.

\smallskip

The second topic of discussion  is a regularity property for all unitary representations of $G$ 
which follows from the G\aa rding inequality  and  is  of significance for a more detailed understanding
of the analytic structure of the group representations.

\begin{prop}\label{plrev4.11} Adopt the assumptions of Proposition~$\ref{plrev4.1}$.
Then  there is a representation independent $a>0$ such that 
\begin{equation}
N_m(\varphi)\leq a\,(\|H\varphi\|+\|\varphi\|)
\label{elrev4.20}
\end{equation}
for all $\varphi\in \chi_m$.
\end{prop}
\proof\
Let $H_1=H^\dagger H$.
Then $H_1$ is a strongly elliptic operator with ellipticity constant $\mu_1\geq \mu^2$.
Then for each  $\lambda\in\langle0,\mu\rangle$ there is a  $\nu\geq0$ such that 
\[
\|H\varphi\|^2=\RRe(\varphi, H_1\varphi)\geq \lambda^2\,N_m(\varphi)^2-\nu^2\,\|\varphi\|^2
\]
for all $\varphi\in \chi_m$ by (\ref{elrev4.1}) applied to $H_1$.
The value of $\nu$ is independent of the particular unitary representation.
Thus 
\[
N_m(\varphi)^2\leq \lambda^{-2}\|H\varphi\|^2+\nu^2\|\varphi\|^2\leq (\lambda^{-1}\|H\varphi\|+\nu\|\varphi\|)^2
\]
for all $\varphi\in \chi_m$.
Now set $a=\lambda^{-1}\vee\nu$.
\hfill$\Box$

\bigskip

The regularity property (\ref{elrev4.20})  was first established for Laplacians  and  unitary representations  by Nelson \cite{Nel}, Section~6, by an algebraic calculation.
(A simplified version of Nelson's result is given by \cite{Rob3}, Lemma~1.7.)
It should be emphasized that the property  is not valid for all representations.
For example, if one considers the representation of $\Ri^d$ by left translations  on $L_p(\Ri^d)$ and sets $H=\Delta$, the standard Laplacian operator,
then Calder{\'o}n \cite{Cal} has shown that  (\ref{elrev4.20}) is valid for all $p\in\langle1,\infty\rangle$.
Nevertheless it fails if $p=1$ or $p=\infty$.
There are  locally integrable functions $\xi$ such that $\Delta\xi$  is also locally integrable but the mixed derivatives $\partial_i\partial_j\xi$ are not.
This pathology has a long history going back at least to  Petrini's 1908 paper \cite{Petr}.
It was, however, still a topical problem in the 60s and the $L_1$ and $L_\infty$ counterexamples can be found in \cite{Orn} and \cite{LM}, respectively.
The situation is even more complicated. 
The Euclidean group is also represented by left translations on the spaces $C(\Ri^d)$ and, more generally, $C^k(\Ri^d)$. 
But  (\ref{elrev4.20}) fails on $C^k(\Ri^d)$ for all integers $k\geq 1$. 

\smallskip
The third topic we address  is a characterization 
 of the analytic elements for general group representations following  suggestions of Roe Goodman  \cite{Goo10} \cite{Goo11}.
 This characterization is in terms of fractional powers of the strongly elliptic operators $H$ and its proof depends
 on the regularity property although the conclusion is independent of the property.
 We next give a brief description of this result.

\smallskip

The earlier discussion of the analytic elements  compared two series with general terms $N_k(\xi)/k!$ and $\|H^k\xi\|/k!$, respectively.
Convergence of the first series characterized the analytic elements $\chi_a(U)$ of the representation $U$ and convergence of the second characterized the analytic elements $\chi_a(H)$ of the strongly elliptic operator $H$.
The arguments of Section~3 established that $\chi_a(H)\subseteq \chi_a(U)$ and consequently $\chi_a(U)$ is dense in the representation space $U$.
But $N_k(\xi)$ involves $k$-derivatives whilst $\|H^k\xi\|$ involves $km$-derivatives.
As Goodman pointed out it is more  appropriate to compare the series with terms  $N_k(\xi)/k!$ and $\|H^k\xi\|/(km)!$.
The latter series is, however, related to the series characterizing the analytic elements of the fractional power $H^{1/m}$ of $H$.
The general theory of  fractional powers of semigroup generators was developed in the late 1950s 
and a summary of  the basic properties  can be found in Chapter~IX of Yosida's book \cite{Yos} on functional analysis or Chapter~1 of Triebel's book on interpolation theory \cite{Tri}.
For current  purposes it suffices to know that if $H$ generates a uniform bounded semigroup then the  fractional powers  $H^\gamma$ with $\gamma\in\langle0,1\rangle$ are well-defined and generate uniformly bounded holomorphic semigroups.
But by the arguments of Section~\ref{S3} each  strongly  elliptic operator $H$ corresponding to a group representation generates 
a continuous semigroup $S$ satisfying operator bounds $\|S_t\|\leq a\,e^{\nu t}$ for some $\nu\geq0$ and all $t>0$.
Therefore the uniform boundedness property can be arranged by replacing $H$ with $H+\nu I$.
This replacement does not change the space of analytic elements of $H$.
Hence  in the following discussion we will assume that $H^{1/m}$ is well-defined and satisfies the standard properties of fractional powers, e.g.\ $(H^{1/m})^k=(H^k)^{1/m}=H^{k/m}$.
Then the analytic elements $\chi_a(H^{1/m})$ of $H^{1/m}$ are defined
 as the $\xi\in \chi_\infty$ such that 
\[
\chi_a(H^{1/m})=\{\xi\in\chi_\infty:\sum_{k\geq 1}\|H^{k/m}\xi\|/k!<\infty\}
\;.
\]
In fact one does not need to consider fractional powers to analyze this subspace.

\begin{lemma}\label{llrev4.1}
The following conditions are equivalent:
\begin{tabel}
\item\label{llrev4.1-1}
$\xi\in \chi_a(H^{1/m})$,
\item\label{llrev4.1-2}
$\sum_{k\geq 1}\|H^k\xi\|/(km)!<\infty$.
\end{tabel}
\end{lemma}
\proof\ 
\ref{llrev4.1-1}$\Rightarrow$\ref{llrev4.1-2}$\;$
If the series defining $\chi_a(H^{1/m})$ is finite then  Condition~\ref{llrev4.1-2} is evident.

\smallskip

\noindent\ref{llrev4.1-2}$\Rightarrow$\ref{llrev4.1-1}$\;$
Since $H^{1/m}$ generates a continuous semigroup there is a $C>0$ such that 
\[
\|H^{l/m}\xi\|\leq C\,(\|H\xi\|+\|\xi\|)
\] 
for all $\xi\in D(H)$ and all $l\in\{1,\ldots, m-1\}$.
Therefore 
\[
\sum_{n\geq 1}\|H^{k/m}\xi\|/k!\leq m\,C\sum_{k\geq0}\Big(\|H^{k+1}\xi\|+\|H^k\xi\|\Big)\Big/(km)!<\infty
\]
for all $\xi$ satisfying Condition~\ref{llrev4.1-2}.
\hfill$\Box$

\bigskip

The final conclusion of Goodman's observations on fractional powers is a complete characterization of the analytic elements
for an arbitrary group representation.

\begin{thm}\label{tlrev4.2}
If $H$ is a strongly elliptic operator associated with the Banach space representation $(\chi, U)$ then $\chi_a(U)=\chi_a(H^{1/m})$.
\end{thm}

This result was established in several  stages.

\smallskip

First Goodman  \cite{Goo10} established the characterization for all unitary representations and  $H$ a Laplacian by a modification
of  Nelson's theory of operator dominance \cite{Nel}.

Secondly, Goodman's  paper also  had a brief but  intriguing  appendix contributed by Nelson that gave an elegant  argument indicating  a similar characterization was valid for  all Banach space representations of the group satisfying the  regularity condition
(\ref{elrev4.20}).
The conclusion of  Nelson's suggestion  was stated for Laplacians in Corollary~A.1 of \cite{Goo10} although Nelson
remarked that his arguments applied to higher orders.

Thirdly, I extended the Goodman--Nelson result to general representations by a longish detour through the interpolation spaces
between the $C^k$-subspaces  $\chi_k$ of the representation space.
My first foray in this direction \cite{Rob1} was again for Laplacians.  
I was aware of Langlands' thesis at that time but was not ambitious enough to extend the interpolation arguments to higher order
operators.
This final step was described in Chapter~II of my book \cite{Robm}.
The reason behind the use of the interpolation arguments was quite simple.
The analytic properties were not affected by transferring to the interpolation spaces but the regularity properties were improved.

The conclusion of the first two stages are summarized  by the following proposition.

\begin{prop}\label{plrev4.2}{\rm (Goodman--Nelson)}
Let $H$ be a strongly elliptic operator associated with the Banach space representation $(\chi, U)$ satisfying 
 $N_m(\xi)\leq a\,(\|H\xi\|+\|\xi\|)$ for some $a>0$ and all $\xi\in \chi_m$.
It follows that  $\chi_a(U)=\chi_a(H^{1/m})$.
\end{prop}

Note that the proposition applies to unitary representations because  unitarity implies the regularity assumption  by
Proposition~\ref{plrev4.11}.
Goodman's proof for unitary representations was based on Nelson's original theory of operator domination \cite{Nel} but the argument
 advanced by Nelson for representations satisfying the regularity assumption was an extension  of this  theory.  
The inclusion $\chi_a(U)\subseteq  \chi_a(H^{1/m})$ is quite elementary and does not depend on domination theory.
For example, if  $\xi\in \chi_\infty$ then there is a $C>0$ such that $\|H^k\xi\|\leq C^k\,\|\xi\|_{km}$ for all $k\geq 1$.
Therefore if $\xi\in \chi_a(U)$ it follows from Lemma~\ref{llrev4.1} that $\xi\in \chi_a(H^{1/m})$.
The proof of the converse inclusion is, however, more delicate.
Nelson's method  was based  on a recursive argument involving the structure relations of the Lie algebra somewhat similar to the reasoning used in Subsection~\ref{S3.3}.
The argument  depends critically  on the regularity condition.
This allows one to estimate products $A^\beta$ with $|\beta|=m$ by  a single $H$, e.g.\ 
if $A^\alpha=A^\beta A^\gamma$ with $|\beta|=m$ then $\|A^\alpha\xi\|\leq C\,(\|HA^\gamma\xi\| +\|A^\gamma\xi\|)$.
Then one can commute  the factor $H$ to the right of the $A^\gamma$ and the additional commutator term 
 is a sum of products $A^\delta$
with $|\delta|\leq |\alpha|-1$, i.e.\ it is a lower order correction.
This is the start of the recursive argument.
In the simplest case,  $G=\Ri^d$, all the $A$ commute and  the regularity condition $\|\xi\|_m\leq C\,(\|H\xi\|+\|\xi\|)$
iterates in this manner to give $\|\xi\|_{km}\leq (2C)^k(\|H^k\xi\|+\|\xi\|)$.
Hence if  $ \xi\in \chi_a(H^{1/m})$ then $\sum_{k\geq1}\|\xi\|_{km}/(km)!<\infty$ by Lemma~\ref{llrev4.1} and this suffices to deduce that
$\xi\in \chi_a(U)$.
Details of the general case are  more complicated since one has to control the lower order terms which arise from the lack of commutation. 
Details are given in Chapter~II of \cite{Robm}. 
We will not persevere with the argument  but instead explain how to deduce Theorem~\ref{tlrev4.2} from Proposition~\ref{plrev4.2}.

\smallskip

First, it  is not surprising that the regularity condition is not necessary.
The Goodman--Nelson arguments essentially use this condition to make a term by term comparison of  the exponential series characterizing  $ \chi_a(H^{1/m})$ with  the series characterizing $ \chi_a(U)$. 
But such a comparison is clearly stronger than necessary for the inclusion  $ \chi_a(H^{1/m}) \subseteq  \chi_a(U)$.
I was aware of this problem by the early 1970s but only realized how to solve it some 10--15 years later.
In the meantime my interests were directed to quite different topics.
My  idea  in the mid 1980s was to exploit the theory of interpolation spaces
and reduce the problem to  a similar problem for an auxiliary representation on a Banach space intermediate to the $C^k$-subspaces $\chi_k$.
In fact it suffices to consider a space intermediate to $\chi$ and  $\chi_1$.

The second observation is that each of the subspaces $\chi_k$ is invariant under the representation $U$.
Hence $U_k=U|_{\chi_k}$ is a representation of $G$ on $\chi_k$. 
Moreover the $C^l$-subspaces $\chi_{k;l}$ of $U_k$ are equal to the $C^{k+l}$-subspaces of $U$, i.e.\ 
 $\chi_{k;l}=\chi_{k+l}$.
Next by a standard procedure of real interpolation (see, for example, \cite{Robm} Section~II.4.1) one can introduce a family of  Banach  spaces
$\chi_{\gamma}$, with $\gamma\in\langle0,1\rangle$, such that  each space is invariant under $U$ and $\chi_1\subseteq \chi_\gamma\subseteq \chi$.
Let   $U_\gamma$ denote the corresponding representations.
The embeddings are continuous and one has bounds $c\,\|\xi\|\leq \|\xi\|_\gamma\leq C\,\|\xi\|_1$ for some $c,C>0$ and all $\xi\in \chi_1$.
Consequently, the $C^k$-subspaces of the representation $U_\gamma$ satisfy $c\,\|\xi\|_k\leq \|\xi\|_{\gamma;k}\leq C\,\|\xi\|_{k+1}$ for all $k\geq 1$ and all $\xi\in \chi_\infty$.
Therefore $\chi_a(U)=\chi_a(U_\gamma)$.
Although there is still a term by term comparison of the two relevant  series there is a slippage of one term in the comparison
which does not affect the conclusion.

Now consider the comparison of the powers of $H$.
For simplicity we use $H$ as a common notation for the operators associated with each of the representations $U$, $U_\gamma$ and $U_1$.
Then one has  $c\,\|H^k\xi\|\leq \|H^k\xi\|_{\gamma}\leq C\,\|H^k\xi\|_1$  for all $\xi\in \chi_\infty$.
But  since $\|\xi\|_1\leq a\,(\|H\xi\|+\|\xi\|)$ it follows  that $\|H^k\xi\|_1\leq a\,(\|H^{k+1}\xi\|+\|H^k\xi\|)$.
Hence the series $\|H^k\xi\|/(km)!$ and $\|H^{k}\xi\|_\gamma/(km)!$ are simultaneously convergent.
Therefore $\chi_a(H^{1/m})=\chi_{\gamma;a}(H^{1/m})$.
Combining these conclusions one obtains the reduction result.

\begin{lemma}\label{llrev4.3}
$\chi_a(U)=\chi_a(H^{1/m})$ if and only if $\chi_{\gamma;a}(U_\gamma)=\chi_{\gamma;a}(H^{1/m})$ for some $\gamma\in\langle0,1\rangle$.
\end{lemma}

It might appear that this manipulation with the interpolation spaces has achieved very little. 
The problem for the representation $(\chi, U)$ has been identified with the analogous problem for the representation $(\chi_\gamma ,U_\gamma)$.
But the redeeming feature, the magic of the interpolation argument,  is that the latter representation satisfies the regularity condition
necessary for  the Goodman--Nelson result, Proposition~\ref{plrev4.2}.
Explicitly, there is an $a_\gamma>0$ such that
\begin{equation}
\|\xi\|_{\gamma;m}\leq a_\gamma\,(\|H\xi\|_\gamma+\|\xi\|_\gamma)
\label{elrev4.3}
\end{equation}
for all $\xi\in \chi_{\gamma;m}$.
Therefore $\chi_{\gamma;a}(U_\gamma)=\chi_{\gamma;a}(H^{1/m})$ by  Proposition~\ref{plrev4.2} and consequently
$\chi_a(U)=\chi_a(H^{1/m})$ by Lemma~\ref{llrev4.3}.
Thus the statement of Theorem~\ref{tlrev4.2} is established.
The only problem remaining is to explain why the regularity property (\ref{elrev4.3}) is valid for the intermediate representations $(\chi_\gamma, U_\gamma)$ even if it is not valid for the representations $(\chi,U)$ and $(\chi_1,U_1)$.
This is a convoluted story.

\smallskip

Interpolation has a long history starting with the work of Riesz in 1926.
But in the late 1950s there was an explosion of interest in the subject motivated by problems
of partial differential operators, approximation theory and singular integration.
Many of the new developments concerned the classical spaces of functions over $\Ri^d$ 
but there were also new ideas on abstract methods of interpolation. 
One of the main motivations for  the construction of new function spaces was indeed the regularity condition (\ref{elrev4.3})  (see \cite{Tri5} pages 38--40).
This led to the construction of various families of spaces satisfying the regularity condition intermediate to the
$C^k$-spaces $L_{p;k}(\Ri^d)$ associated with left translations on the $L_p$-spaces.
So the proof of Theorem~\ref{tlrev4.2} for the $\Ri^d$-theory could have been completed by appealing to the results described,
for example, in the books of Triebel \cite{Tri} \cite{Tri5}.
Unfortunately, there was no equivalent theory for representations of Lie groups, although Peetre gave some 
partial results in \cite{Pee1}.
Nevertheless the methods required to describe the general situation were available.
In particular there was a detailed description of interpolation methods for semigroups acting on abstract Banach spaces
in the book by Butzer and Berens \cite{BB} on approximation theory.
This theory was largely based on ideas of Peetre on methods of real interpolation between general Banach spaces.
We now sketch its application to the representation $(\chi, U)$ of the Lie group $G$.

First, if $\chi_1$ is   the $C^1$-subspace of the representation space $\chi$ then
$\chi_{\gamma,p}=(\chi,\chi_1)_{\gamma, p}$ is defined as the space of $\xi\in \chi$ such that the seminorm $N_{\gamma,p}(\xi)=(\int^\infty_0 dt\, t^{-1}\,(t^{-\gamma}\kappa_t(\xi))^p)^{1/p} $ is finite where $\kappa_t(\xi)=\inf_{\xi_1\in \chi_1}(\|\xi-\xi_1\|+t\,\|\xi_1\|_1)$ and $p\in[1,\infty\rangle$.
The interpolation function $\kappa_t$  gauges  the relevant importance of the representation space  $\chi$ and the $C^1$-subspace $\chi_1$.
Consequently $\gamma$ gives a measure of the smoothness of $\xi$ with the choice of $p$ giving an extra gradation.
If $\xi\in\chi_1$ then $\kappa_t(\xi)$ tends to zero as  $t\to 0$.
Secondly,   if $S$ is the  continuous semigroup generated by the strongly elliptic operator $H$ associated with $(\chi,U)$
then $\chi^S_{\gamma, p}$  is defined  as the subspace of  $\xi\in\chi$ for which  the seminorm $N^S_{\gamma,p}(\xi)=(\int^\infty_0 dt \,t^{-1}\,( t^{-\gamma}\|(I-S_t)\xi\|)^p)^{1/p}$ is finite.
Since  the semigroup $S$ is holomorphic 
it also  follows that $\chi^S_{\gamma, p}$ is the subspace of $\xi$ for which the seminorm $(\int^\infty_0\, dt \,t^{-1}\,(t^{1-\gamma}\|HS_t\xi\|)^p)^{1/p}$ is finite.
In these definitions the $\gamma$ and $p$ measure the smoothing properties of the semigroup $S$ for small~$t$.
The only apparent group connection between  $\chi_{\gamma,p}$ and $\chi^S_{\gamma,p}$ is the first space involves  the $C^1$-subspace $\chi_1$ of the representation space whilst the second  depends indirectly on the representation $U$ through the strongly elliptic operator $H$.
Nevertheless the  two intermediate spaces both give a measure of smoothness and the 
striking conclusion is  that  they are equal, modulo a slight change of parameter.
Specifically,  $\chi_{\gamma, p}=\chi^S_{(\gamma/m), p}$ with equivalence of the natural norms.
In fact there is even a third chacterization of these spaces directly involving the representation $U$.
The space $\chi_{\gamma,p}$ consists of the $\xi\in\chi$ for which the seminorm $\int_\co dg\,|g|^{-d}(|g|^{-\gamma}\|(I-U(g)\xi\|)^p$
is finite where $\co$ is an open neighbourhood of the identity in $G$.
These results are given by Proposition~II.4.3 and Theorem~II.6.1 of \cite{Robm} although they were well known for $G=\Ri^d$ much earlier.

Once one has the identification  $\chi_{\gamma, p}=\chi^S_{(\gamma/m), p}$ it is relatively straightforward to deduce the regularity property (\ref{elrev4.3}) for the intermediate spaces $\chi_\gamma=\chi_{\gamma, p}$.
One key observation is that $\xi_1=S_t\xi\in\chi_1$ for each $\xi\in\chi$.
Therefore the decomposition $\xi=(\xi-\xi_1)+\xi_1$ takes the form $\xi=(I-S_t)\xi+S_t\xi$ and allows one to estimate $\kappa_t(\xi)$ in terms
of $\|(I-S_t)\xi\|$ and $\|S_t\xi\|_1$.
The details are given in Theorem~II.4.5 of \cite{Robm} but the ideas are just borrowed from the $\Ri^d$-theory developed in the 1960s
which can be found in \cite{BB} or \cite{Tri} among many other places.

\smallskip

This completes our discussion of  the characterization of analytic elements, Theorem~\ref{tlrev4.2}, and our summary  of the developments concerning the higher order strongly elliptic operators introduced by Langlands.
The most striking aspect of these results is their universal nature, e.g.\  $\chi_a(U)=\chi_a(H^{1/m})$ for all the suitably normalized  $m$-th order operators $H$  independent of the group structure.
The conclusions are basically locally and are independent of the Lie algebraic details.
One can obtain  more detailed global results by specializing to second-order operators such as Laplacians or to restricted classes 
of groups.
But the global analysis  requires the introduction of quite different techniques, e.g.\ generalized Nash inequalities \cite{Robm} or Harnack inequalities \cite{VSC}.
The conclusions  are also sensitive to the large scale geometry of the group.
A detailed analysis of these properties for groups of polynomial growth can be found in \cite{DER4}.

We conclude with a brief discussion of a slightly different topic, general order subelliptic operators.

The  subelliptic theory is formulated in a similar manner to the strongly elliptic theory but the vector space basis $a_1,\ldots, a_d$ of the Lie algebra $\gotg$ is replaced by an algebraic subbasis $a_1, \ldots, a_{d_1}$, i.e.\ a linearly independent set of elements which generate $\gotg$ algebraically.
Then one can define differential operators as polynomials of  the representatives $A_k$ of the $a_k$ in the subbasis.
The properties of second-order subelliptic operators, so-called `sums of squares', have been extensively studied since H\"ormander's fundamental 1967 paper \cite{Hor1}.
If, however, one tries to develop the structure of  higher order subelliptic operators following the outlines of the strongly elliptic theory one immediately encounters several new obstacles.

The first  substantial obstacle is to find  a  replacement for the strong ellipticity condition (\ref{elrev2.0}).
This is an $\Ri^d$-condition on the coefficients which does not reflect the  restraints imposed by the subellipticity condition.
Secondly, it is not at all clear that there is an alternative version of the parametrix arguments.
This  problem is related to the previous difficulty.
The parametrix expansion for operators associated with the group $G$ is  the  analogue of perturbation theory with  the  unperturbed  system given by operators corresponding to $\Ri^d$.
Fortunately both these obstacles can be avoided by a rather different common approach.
The basic idea is to introduce a `simpler' group $G_0$, related to $G$ but with a streamlined algebraic structure dictated by the subelliptic basis, as a replacement for $\Ri^d$.

First, however, define the subelliptic distance $|g|_1$ as  the shortest length  of the absolutely continuous paths from $g$ to $e$ following the directions of the algebraic subbasis.
It is not evident that one can find connecting paths of this type for each $g\in G$  but this is a result of Caratheodory's early research  into thermodynamics \cite{Car}.
Moreover, it follows that a ball of radius $\delta$ measured with respect to this distance behaves as $\delta^D$ as $\delta\to0$ where 
$D$ is an  integer, the local subelliptic dimension.
It can be calculated as follows.
Let ${\gotg}_1$ denote the linear span of the algebraic basis 
$a_1,\ldots ,a_{d_1}$ and ${\gotg}_j$ the span of the algebraic basis together
with the corresponding multiple commutators of order less
than or equal to $j$. 
Then ${\gotg}_1 \subset
{\gotg}_2 \subset \ldots \subset {\gotg}_{r} = {\gotg}$ where $r$ is an integer, the rank of the algebraic basis.
Next set $\gotg_1' = {\gotg}_1$, and 
$\gotg_j '$ the vector space  complement of ${\gotg}_{j-1}$ in  ${\gotg}_j$.  
This yields the direct sum decomposition
${\gotg} =\gotg_1' \oplus \gotg_2'\oplus
\cdots \oplus\gotg_r'$ of the Lie algebra.
Then $D$ is given by $D=\sum^r_{j=1} j\,(\dim{\gotg_j'})$.

Secondly, 
the group $G_0$ is is defined by a contraction procedure.
Define $\gamma$ as the family of maps of $\gotg$ into $\gotg$ such that  $\gamma_t(a)=t^k a$ for all $a=\gotg_k'$ and $t>0$.
Then $\gotg_0$ is defined as the vector space $\gotg$ equipped with the Lie bracket
\[
[a,b]_0=\lim_{t \to0}\gamma^{-1}_t([\gamma_t(a),\gamma_t(b)])
\;.
\]
It follows that $\gotg_0$ is a nilpotent Lie algebra and 
$a_1,\ldots,a_{d_1}$   is an algebraic basis of $\gotg_0$ of  rank $r$.
Moreover, the dilations $\gamma_t$ are automorphisms of $\gotg_0$.
Then  $G_0$ is defined as the connected, simply connected, Lie group with Lie algebra $\gotg_0$.
It is this group  which acts as the local approximation to $G$ in the subelliptic
theory.
The simplifying feature of $G_0$ is the existence of the dilations $\gamma_t$
which allow scaling arguments to extend local properties globally.
Note that  if $a_1,\ldots, a_{d_1}$ is a vector space basis of
$\gotg$ then $\gotg_1=\gotg$ and $\gamma_t(a)=ta$ for all $a\in\gotg$.
Therefore $\gamma^{-1}_t([\gamma_t(a),\gamma_t(b)])=t\,[a,b]\to0$ as $t\to0$ and $\gotg_0$ is
abelian. Thus $G_0=\Ri^d$ in conformity with the earlier  strongly elliptic case.

Thirdly, the notion of subellipticity of the operator $H=\sum_{\alpha: |\alpha|\leq m} c_\alpha A^\alpha$, where the multi-indices $\alpha$
only involve the indices $\{1,\ldots, d_1\}$ of the subbasis, is defined in a manner that simulates the definition of strong ellipticity.
Since Corollary~\ref{clrev4.10} establishes that strong ellipticity is equivalent to the G\aa rding inequality (\ref{elrev4.1}) for 
left translations on $L_2(\Ri^d)$ we define the operator $H$ to be subelliptic on $G$ if  an analogous inequality is satisfied on 
$C_c(G_0)$.
Explicitly, $H$ is subelliptic on $G$  if 
\begin{equation}
\RRe(\varphi, H\varphi)\geq \lambda\, N'_{m/2}(\varphi)^2-\nu\,\|\varphi\|^2
\label{elrev4.11}
\end{equation}
for all $\varphi\in C_c^\infty(G_0)$ where $N'_k$ is the seminorm given by restricting the supremum in the earlier definition of $N_k$ to multi-indices in the subelliptic directions $a_1,\ldots, a_{d_1}$.
Therefore  $H$ is subelliptic on $G$ if and only if it is subelliptic on $G_0$.
This definition removes the first obstacle cited above.

Next, as a preliminary to developing a parametrix formalism for the subelliptic operators in general representations of $G$,  one must first analyze the operators in the left regular representation of $G_0$ on  $L_2(G_0)$.
In the strongly elliptic case with $G_0=\Ri^d$ this was accomplished in Section~\ref{S2}  largely by techniques of Fourier analysis.
In the subcoercive case the situation is more complicated.
It is, however, facilitated by  the nilpotent structure and the homogeneity properties with respect to dilations.
In particular one establishes that each closed  subelliptic operator $H$ generates a continuous semigroup $S$ on the $L_p$-spaces over $G_0$ with a kernel satisfying $m$-th order Gaussian bounds similar to those given by (\ref{elrev2.910})
but with $|g|$ replaced by $|g|_1$ and $d$ replaced by $D$.
Subsequently the properties of $H$ for the nilpotent group $G_0$ are extended to the corresponding operator on the group $G$ by  the parametrix arguments.
The reasoning is not substantially different.
Although  the  conclusions for the  semigroup structure in the subelliptic case are  directly analogous to those of the strongly elliptic case their implications for the differential and analytic structure are considerably weaker.
One striking difference is the failure of the Goodman characterization of the analytic functions in terms of fractional powers of the 
elliptic operators.

If the semigroup $S$ generated by the $m$-th order subelliptic operator $H$ is uniformly bounded one can define $H^{1/m}$ as 
before and the characterization of the subspace of  analytic elements $\chi_a(H^{1/m})$ given by Lemma~\ref{llrev4.1} is still valid.
Hence $\chi_a'(U)\subseteq \chi_a(H^{1/m})$ where  $\chi_a'(U)$ is the subspace of analytic elements of the representation $U$ defined with the subelliptic seminorms $N_k'$.
It is not, however, true that  $\chi_a'(U)=\chi_a(H^{1/m})$ even for second-order operators and unitary representations.
Example~8.7 of \cite{ER1} gives a counterexample based on the left regular representation of the group of rotations on $\Ri^3$
and the standard Laplacian.
Nevertheless  many regularity results have been established in the subelliptic case.
Details can be found in \cite{ER7} \cite{ER8} and \cite{ERS}.

Finally we note that there is a third class of elliptic operators which can be analyzed by Langlands' methods, weighted strongly elliptic operators.
Subellipticity is based on the idea that there is a certain subset of preferred directions. 
In the weighted theory all directions are allowed but some have greater weight, or preference, than others.
An extensive  analysis of this class of operators, along the foregoing lines, can be found in \cite{ER9}.
In particular there is an analogous theory of holomorphic semigroups generated by weighted operators.
Again this leads to a good understanding of the corresponding differential structure and 
the structure of the weighted analytic elements.
For example  the characterization of the analytic elements in terms of the corresponding elements of fractional powers remains valid.
Nevertheless there are significant differences introduced by the weighting.
In conclusion the structural properties   of the analytic elements  in the broader context of subelliptic operators or weighted elliptic  operators  still pose intriguing open  problems 60 years after Langlands' thesis work.

\bigskip

\bigskip
\noindent
\begin{tabular}{@{}cl@{\hspace{10mm}}cl}
&   Derek W. Robinson   &  {} &{}\\
&   Mathematical Sciences Institute (CMA)    &  {} &{}\\
  &Australian National University& & {}\\
&Canberra, ACT 0200 && {} \\
  & Australia && {} \\
  &derek.robinson@anu.edu.au
 & &{}\\
\end{tabular}

\end{document}